\documentclass[11pt, leqno]{article}
\usepackage{amsmath,amsfonts,amssymb,wasysym}%,hyperref} %amsthm
\usepackage[latin1]{inputenc}
\usepackage{graphicx}
\usepackage{color}
\newcommand{\ren}{{\mathbb R}^N}

\newcommand{\be}[1]{\begin{equation}\label{#1}}
\newcommand{\ee}{\end{equation}}

\newcommand{\prf}{\par\smallskip\noindent{\sl Proof. \/}}
\newcommand{\finprf}{\unskip\null\hfill$\;\square$\vskip 0.3cm}

\newtheorem{theorem}{Theorem}[section]
\newtheorem{lemma}{Lemma}[section]

\newtheorem{proposition}[theorem]{Proposition}

\newcommand{\ve}{\varepsilon}
\numberwithin{equation}{section}
 
\def\qed{\,\unskip\kern 6pt \penalty 500
\raise -2pt\hbox{\vrule \vbox to8pt{\hrule width 6pt
\vfill\hrule}\vrule}\par}
\definecolor{darkblue}{rgb}{0.05, .05, .65}
\definecolor{darkgreen}{rgb}{0.1, .65, .1}
\definecolor{darkred}{rgb}{0.8,0,0}
\newcommand{\RN}{\mathbb{R}^N}
\newcommand{\RR}{\mathbb{R}}
\numberwithin{equation}{section}

\definecolor{darkblue}{rgb}{0.05, .05, .65}
\definecolor{darkgreen}{rgb}{0.05, .55, .05}
\definecolor{darkred}{rgb}{0.8,0,0}

\begin{document}

\title{\bf The Fisher-KPP equation \\ with nonlinear fractional diffusion}

\author{Diana Stan$^{\,*}$
~and~ Juan Luis V\'azquez \footnote{ Departamento de Matem\'aticas, Universidad Aut\'onoma de Madrid, 28049 Madrid, Spain. \newline   E-mails: {\tt
diana.stan@uam.es,  juanluis.vazquez@uam.es}
}
 \\}

\maketitle

\begin{abstract}
We study the propagation properties of nonnegative and bounded solutions of the class of reaction-diffusion equations with nonlinear fractional diffusion: $u_{t} + (-\Delta)^s (u^m)=f(u)$. For all $0<s<1$ and $m> m_c=(N-2s)_+/N $, we consider the solution of the initial-value problem  with
initial data having fast decay at infinity  and prove that its level sets propagate exponentially fast in time, in contradiction to the traveling wave behaviour of the standard KPP case, which corresponds to
putting $s=1$, $m=1$ and $f(u)=u(1-u)$. The proof of this fact uses as an essential ingredient the recently established decay properties of the self-similar
solutions of the purely diffusive equation,  $u_{t} + (-\Delta)^s u^m=0$.
\end{abstract}

    \vspace{1cm}


%35A05, %General existence and uniqueness theorems
%35K55, %Nonlinear PDE of parabolic type
%35S10, %Initial value problems for PsDO

\medskip


%%%%%%%%%%%%%%%%%%%%%%%%%%%%%%%%%%%%%%%%%%%%%%%%%%%%%%%%%%%%%%%%%%%%%%%%%%%%%%

\newpage
\tableofcontents

\newpage

\section{Introduction}\label{sect.intro}

We consider the following reaction-diffusion problem
\begin{equation}\label{KPP}
  \left\{ \begin{array}{ll}
  u_{t}(x,t) + L_s u^m(x,t)=f(u) &\text{for } x \in \mathbb{R}^N \text{ and }t>0, \\
  u(x,0)  =u_0(x) &\text{for } x \in \mathbb{R}^N,
    \end{array}
    \right.
\end{equation}
where $L_s=(-\Delta)^s$ is the Fractional Laplacian operator with $s \in (0,1)$. We are interested in studying the propagation properties of nonnegative and
bounded solutions of this problem in the spirit of the Fisher-KPP theory. Therefore, we assume that the reaction term $f(u)$ satisfies
\begin{equation}\label{propf}
f\in C^1([0,1])\text{ is a concave function with } f(0)=f(1)=0 , \quad f'(1)<0<f'(0).
\end{equation}
For example we can take $f(u)=u(1-u).$  Our results will depend on the parameters $m $ and $s$, according to the ranges $m_c<m<m_1$, $m_1<m\le 1$, and $m>1$,
where
$$
 m_c=\frac{(N-2s)_+}{N}, \qquad m_1=\frac{N}{N+2s}\,.
$$

\noindent {\bf 1.1. Perspective. The traveling wave behavior. } \noindent The problem with standard diffusion goes back to the work of Kolmogorov, Petrovskii and Piskunov, see \cite{KPP}, that presents the most simple reaction-diffusion equation concerning the concentration $u$ of a single substance in one spatial dimension,
\begin{equation}\label{KPPdim1}
\partial_t u=D u_{xx} + f(u).
\end{equation}
 The choice $f(u) = u(1-u)$ yields Fisher's equation \cite{Fisher} that was originally used to describe the spreading of biological populations. The celebrated
 result says that the long-time behavior of any solution of \eqref{KPPdim1}, with suitable data $0\le u_0(x)\le 1$ that decay fast at infinity, resembles a traveling
 wave with a definite speed. When considering equation \eqref{KPPdim1} in dimensions $N\geq 1$, the problem becomes
\begin{equation}\label{classicalKPP}
u_t-\Delta u=f(u) \quad \text{in }(0,+\infty)\times \RN,
\end{equation}
which corresponds to \eqref{KPP} in the case when $L_s=-\Delta$,  the standard Laplacian. This case has been studied by Aronson and Weinberger in
\cite{AronsonWeinberger2,AronsonWeinberger}, where they prove the following result.

%%%%%%%%%%%%%%%%%%%%%%
\noindent  {\bf Theorem AW.} {\sl Let $u$ be a solution of \eqref{classicalKPP} with $u_0 \neq 0$ compactly supported in $\RN$ and satisfying $0\leq u_0(\cdot)\leq 1$. Let
$c_*=2\sqrt{f'(0)}$. Then,
\begin{enumerate}
\item if $c>c_*$, then $u(x,t) \rightarrow 0$ uniformly in $\{|x|\geq ct\}$ as $t\rightarrow \infty$.
\item if $c<c_*$, then $u(x,t) \rightarrow 1$ uniformly in $\{|x|\leq ct\}$ as $t\rightarrow \infty$.
\end{enumerate}
}

%%%%%%%%%%%%%%%%%%%%%%

\noindent In addition, problem \eqref{classicalKPP} admits planar traveling wave solutions connecting $0$ and $1$, that is, solutions of the form
$u(x,t)=\phi(x\cdot e+ct)$ with
$$
-\phi'' + c\phi'=f(\phi) \text{ in } \RR, \quad \phi(-\infty)=0, \ \phi(+\infty)=1.
$$

This asymptotic traveling-wave behavior  has been generalized in many interesting ways. Of concern here is the consideration of nonlinear diffusion. De Pablo and
V\'{a}zquez study in \cite{dPJLVjde1991} the existence of traveling wave solutions and the property of finite propagation for the reaction-diffusion equation
$$u_t=(u^m)_{xx}+\lambda u^n(1-u), \quad (x,t) \in \mathbb{R}\times (0,\infty)
$$
with $m>1$, $\lambda>0$, $n \in \RR$ and $u=u(x,t)\geq 0$. Similar results hold also for other slow diffusion cases, $m>1$, studied by de Pablo and S\'{a}nchez
(\cite{dePabloSanchez98}).

\medskip

\noindent {\bf \bf 1.2. Non-traveling wave behavior.} \noindent Departing from these results, King and McCabe  examined in \cite{KingMcCabe} a case of fast diffusion,
namely
$$
u_t=\Delta u^{m}  +u(1-u),\quad x \in \mathbb{R}^N, t>0,
$$
where $(N-2)_+/N<m<1$. They showed  that the problem does not admit traveling wave solutions. Using a detailed formal analysis, they also showed that level sets of the
solutions of the initial-value problem with suitable initial data propagate exponentially fast in time. \normalcolor They extended the results to all $0<m<1$.

On the other hand, and independently,  Cabr\'{e} and Roquejoffre in \cite{CabreRoquejoffre1, CabreRoquejoffreArxiv} studied the case of fractional linear diffusion,  $s\in (0,1)$ and  $m=1$, and they concluded in the same vein that there is no traveling wave behavior as $t\to\infty$, and indeed the level sets propagate exponentially fast in time. This came as a surprise since their problem deals  with linear diffusion.

Motivated by these two examples of break of the asymptotic TW structure, we  study here the case of a diffusion that is both fractional and nonlinear, namely
problem \eqref{KPP} in the range $s\in (0,1)$ and  $m>m_c$. The initial datum $u_0(x) : \mathbb{R}^N \rightarrow [0,1]$ and satisfies a growth condition of the
form
\begin{equation}\label{dataAssump}
0\le u_0(x) \leq C|x|^{-\lambda(N,s,m)}, \quad \forall x \in \RN,
\end{equation}
where the exponent $\lambda(N,s,m)$ is stated explicitly in the different ranges,  $m_c<m<m_1$ and $m_1<m$. In this paper we establish the negative result about traveling wave behaviour,  more precisely, we prove that an exponential rate of propagation of level sets is true in all cases. We also explain the mechanism for it in simple terms: the exponential rate of propagation of the level sets of solutions (with initial data having a certain minimum decay for large $|x|$) is a consequence of the power-like decay behaviour of the fundamental solutions of the diffusion problem studied in \cite{VazquezBarenblattFractPME}. Therefore, we obtain two main cases in the analysis,  $m_c<m<m_1$ and $m>m_1$, depending on that behaviour.

\medskip

\noindent {\bf 1.3. Main results.} The existence of a unique mild solution of problem \eqref{KPP} follows by semigroup approach. The mild solution corresponding to an initial datum $u_0 \in
L^1(\RN),$ $0\leq u_0\leq1$ is in fact a positive, bounded, strong solution with $C^{1,\alpha}$ regularity. In the Appendix we give a brief discussion of these
properties. Let us introduce some notations. Once and for all, we put   $\beta=1/(N(m-1)+2s)$ \ and
\begin{equation}\label{sigma}
\sigma_1 =\frac{1-m}{2s}f'(0), \quad  \sigma_{2}=\frac{1}{N+2s}f'(0), \quad \sigma_3=\frac{1+2(m-1)\beta s}{N+2s}f'(0).
\end{equation}
The value $\sigma_1$ appears for $m_c<m<m_1$ and then $\sigma_1>\sigma_2$.
Notice also that  $\sigma_{2}<\sigma_3$ for $m>1$. Here is the precise statement of our main results for the solutions of the generalized KPP problem
\eqref{KPP}.

\begin{theorem}\label{mainThm1}
Let $N\geq 1$, $s\in (0,1)$, $f$ satisfying \eqref{propf} and $m_1<m \leq 1$. Let $u$ be a solution of \eqref{KPP}, where $0\leq u_0(\cdot)\leq 1$ is measurable,
$u_0 \neq 0$ and satisfies
\begin{equation}
0\le u_0(x) \leq C|x|^{-(N+2s)}, \quad \forall x \in \RN.
\end{equation}
 Then
\begin{enumerate}
\item if $\sigma>\sigma_{2}$, then $u(x,t) \rightarrow 0$ uniformly in $\{|x|\geq e^{\sigma t} \}$ as $t\rightarrow \infty$.
\item if $\sigma<\sigma_{2}$, then $u(x,t) \rightarrow 1$ uniformly in $\{|x|\leq e^{\sigma t} \}$ as $t\rightarrow \infty$.
\end{enumerate}
\end{theorem}

\begin{theorem}\label{mainThm2}
Let $N\geq 1$, $s\in (0,1)$, $f$ satisfying \eqref{propf} and $m_s<m<m_1$. Let $u$ be a solution of \eqref{KPP}, where $0\leq u_0(\cdot)\leq 1$ is measurable,
$u_0 \neq 0$ and satisfies
\begin{equation}
0\leq u_0(x) \leq C |x|^{-2s/(1-m)},  \quad \forall x \in \RN.
\end{equation}
  Then
\begin{enumerate}
\item if $\sigma>\sigma_{1}$, then $u(x,t) \rightarrow 0$ uniformly in $\{|x|\geq e^{\sigma t} \}$ as $t\rightarrow \infty$.
\item if $\sigma<\sigma_{1}$, then $u(x,t) \rightarrow 1$ uniformly in $\{|x|\leq e^{\sigma t} \}$ as $t\rightarrow \infty$.
\end{enumerate}
\end{theorem}

\begin{theorem}\label{mainThm3}
Let $N\geq 1$, $s\in (0,1)$, $f$ satisfying \eqref{propf} and $m>1$. Let $u$ be a solution of \eqref{KPP}, where $0\leq u_0(\cdot)\leq 1$ is measurable, $u_0 \neq
0$ and satisfies $$0\le u_0(x) \leq C|x|^{-(N+2s)}, \quad \forall x \in \RN.$$ Then
\begin{enumerate}
\item if $\sigma>\sigma_{3}$, then $u(x,t) \rightarrow 0$ uniformly in $\{|x|\geq e^{\sigma t} \}$ as $t\rightarrow \infty$.
\item if $\sigma<\sigma_{2}$, then $u(x,t) \rightarrow 1$ uniformly in $\{|x|\leq e^{\sigma t} \}$ as $t\rightarrow \infty$.
\end{enumerate}
\end{theorem}

\begin{figure}[h!]\label{figparameters}
  \centering
  \includegraphics[width=80mm,height=60mm]{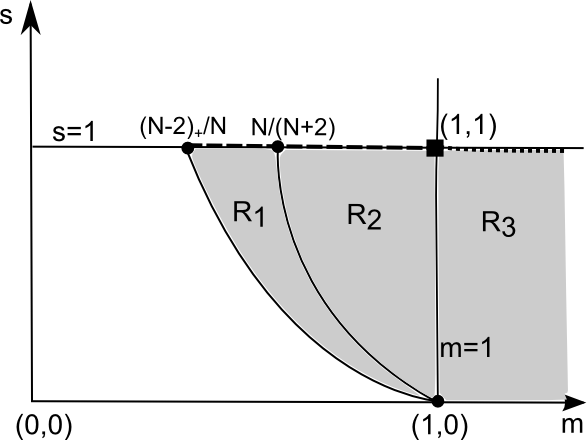}
  \caption{Ranges of parameters $m$ and $s$: we study the cases $R_1=\{s\in(0,1), \ (N-2s)_+/N<m<N/(N+2s)\}$; $R_2 =\{s\in(0,1), \ N/(N+2s)<m\leq 1\} $,
  $R_3 =\{s\in(0,1), \ m> 1\} $}
\end{figure}

\medskip

\noindent \textbf{Remarks.} In all ranges of parameters $m>m_c$, there appear critical values of $\sigma$ with an influence on the behavior of the level sets.

\noindent$\bullet$ In the case $m_1<m<1$, the case $\sigma=\sigma_2$ is still open. This critical exponent is the same as in the case of the linear diffusion $m=1$, proved in \cite{CabreRoquejoffreArxiv}.

\noindent$\bullet$ In the range $m_c<m<m_1$, the case $\sigma=\sigma_1$ is still open. In particular, for the classical case $s=1$ and $f(u)=u(1-u)$ we get $\sigma_1=\frac{1-m}{2}$, which is a critical speed found by  King and McCabe \cite{KingMcCabe}. In this way, we complete their result with rigorous proofs to all $s\in (0,1)$.

\noindent$\bullet$ In the case $m>1$, we do not cover the entire interval $[\sigma_{2},\sigma_{3}]$. Therefore, we prove that for $m>1$ the nonlinearity has a
different influence on the velocity of propagation.

\noindent$\bullet$ The result of Theorems \ref{mainThm1} and \ref{mainThm2} is true also in the case $m=m_1$, where $\sigma_1=\sigma_2$. The outline of the proof is the same, but there are a number of additional technical difficulties, typical of borderline cases. We have decided to skip the lengthy analysis of this case because of the lack of novelty for our intended purpose.

\medskip

\noindent   Our main conclusion is  that  exponential propagation is shown to be the common occurrence, and the existence of traveling wave behavior is reduced to the classical KPP cases mentioned at the beginning of this discussion (see dotted line in Figure \ref{figparameters}).

\medskip

As we have already mentioned, one of the  motivations of the work was to make clear the mechanism that explains the exponential rate of expansion in simple terms, even in this situation that is more complicated than \cite{CabreRoquejoffre1, CabreRoquejoffreArxiv}. In fact, due to the nonlinearity, the solution of the diffusion problems involved in the proofs does not admit an integral representation as the case $m=1$. Instead, we will use as an essential tool the behavior of
the fundamental solution of the Fractional Porous Medium Equation, also called Barenblatt solution, recently studied in \cite{VazquezBarenblattFractPME}. To be
precise, the decay rate of the tail of these solutions as $|x|\to\infty$ is the essential information we use to calculate the rates of expansion. This information
is combined with more or less usual techniques of linearization and comparison with sub- and super-solutions.  We also need accurate lower estimates for positive solutions of this latter equation, and a further selfsimilar analysis for the linear diffusion problem.

\medskip

\noindent {\bf 1.4. Organization of the proofs.}  In Section \ref{SubsecConv0}, under the assumption of initial datum with the decay \eqref{dataAssump}, we prove convergence to $0$ in the outer set $\{|x|\geq e^{\sigma t} \}$ by constructing a super-solution of the linearized problem with reaction term $f'(0)u$. The arguments hold for
$\sigma$ larger than the corresponding critical velocity.

In Section \ref{SubsecConv1} we prove convergence to $1$ on the inner sets $\{|x|\leq e^{\sigma t} \}$ in various steps. We only assume $0\leq u_0\le 1$, $u_0
\neq 0.$ We first show that the solution reaches a certain minimum profile for positive times, thanks to the analysis of Theorem \ref{low.par.est1} below, we then perform an iterative proof the conservation in time of this minimum level, and finally
 convergence to $1$ is obtained by constructing a super-solution to the problem satisfied by $1-u^m$. Therefore, we deal with a problem of the form
$$
a(x,t)\,w_t (x,t) + L_s w(x,t)+b_0 \,w(x,t)\ge 0.
$$
A suitable choice for constructing the super-solution $w$ is represented by self-similar solutions of the form $U(x,t)=t^{\alpha'}F(|x|t^{-\beta'})$ of the linear problem
\begin{equation}\label{LinProblem2}
U_t +L_s U=0
\end{equation}
 with radial increasing initial data. This motivates us to derive a number of properties of the linear diffusion problem \eqref{LinProblem2}, also known as the Fractional Heat Equation. In particular, we need to show that the profile $F$ mentioned above has the same asymptotic behavior as the initial data. In order to establish such fact we have to review,  Section \ref{sectSelfSimLinProb},  the properties of the fundamental solution of Problem \eqref{LinProblem2}
 $$
 K_s(x,t)=t^{-\frac{N}{2s}}f(t^{-\frac{1}{2s}}|x|), \qquad f(r) \sim r^{-(N+2s)}.
 $$
 We perform a further analysis of the profile $f$ by proving that $rf' \sim r^{-(N+2s)}.$

%This analysis also  applies in the case $m_c<m<m_1$, $\sigma<\sigma_1$, by using the corresponding tail behavior $|x|^{-2s/(1-m)}$. This is done in Section
%\ref{subsectFirstRange}.

\medskip

\noindent\textbf{Remark.} As a consequence of the exponential propagation of the level sets, we immediately obtain the non-existence of traveling wave solutions
of the form $u(x,t)=\varphi(x+t\cdot e)$.
However, our results amount to the existence of a kind of logarithmic traveling wave behaviour, that is a kind of wave solutions that travel linearly if we measure distance in a logarithmic scale. This whole issue deserves further investigation.

\medskip

\noindent \textbf{1.5. New estimates for the fractional diffusion problem. } \noindent The study of the sub- and super-solutions is strongly determined by the
existence of suitable lower parabolic estimates for the associated diffusion problem, the Fractional Porous Medium Equation (FPME)
\begin{equation}\label{FPME}
  \left\{ \begin{array}{ll}
  \underline{u}_{t}(x,t) + L_s \underline{u}^m(x,t)=0 &\text{for } x \in \mathbb{R}^N \text{ and }t>0, \\
  \underline{u}(x,0)  =\underline{u}_0(x) &\text{for } x \in \mathbb{R}^N.
    \end{array}
    \right.
\end{equation}
In Section \ref{SubsectLowerParEst}, we devote a separate study in the case $m>1$ of the behavior of the solution when $|x|\rightarrow \infty$, more precisely its
rate of decay, for small times $t>0$.  Our main result says that roughly speaking
\begin{equation*}
\underline{u}(x,t) \sim \,t\,|x|^{-(N+2s)}
\end{equation*}
when $|x|$ is large and $t$ small. The precise result is as follows.

\begin{theorem}\label{low.par.est1}
Let $\underline u(x,t)$ be a solution of Problem \eqref{FPME} with initial data $u_0(x)\ge0$ such that $u_0(x)\ge 1$ in the ball $B_1(0)$. Then there is a time
$t_1>0$ and constants $C_*, R>0$ such that
\begin{equation}
\underline{u}(x,t)\ge C_*\,t\,|x|^{-(N+2s)}
\end{equation}
if $|x|\ge R$  and $0<t<t_1$.
\end{theorem}

The fact that solutions of the FPME with nonnegative initial data become immediately positive for all times
$t>0$ in the whole space has been proved in \cite{DPQRV1,DPQRV2}. Such result is true not only for $0<s<1$ and $m>1$, but also for $0<s<1$ and $m>m_c=(N-2s)_+/N$,
this lower restriction on $m$ aimed at avoiding the possibility of extinction in finite time.

Precise quantitative estimates of positivity for $t>0$ on bounded domains of $\ren$ have been  obtained in the recent paper \cite{BV2012}.
The estimates of that reference are  also precise in describing the behavior as $|x|\to\infty$ when $m<1$ (fast diffusion),
but they are not relevant to establish the far-field behavior for $m>1$. We recall that the limit $s\to1$ with $m>1$ fixed we get the standard porous medium
equation, where positivity at infinity for all nonnegative solutions is false due to the property of finite propagation, cf. \cite{VazPMEwholespace}. This
explains that some special characteristic of fractional diffusion must play a role if positivity is true.

We fill the needed gap  for some convenient class of initial data that includes continuous nontrivial and nonnegative initial data. We give a quantitative version
here since it can be useful in the applications.

%%%%%%%%%%%%%%%%%%%%%%%%%%%%%%%%%%%%%%%%%%%%%%%%%%%%%%%%%%%%%%%%%%%%%%%%%%%%%%%%%%%%%%%

\section{Preliminaries}

\subsection{Nonlinear diffusion. The Fractional Porous Medium Equation}

We recall some useful results concerning the porous medium equation with fractional diffusion (FPME). We refer to \cite{DPQRV2} where the authors develop the
basic theory  for the general problem
\begin{equation}\label{genfractPME}
  \left\{ \begin{array}{ll}
 u_t= -L_s (|u|^{m-1}u) &\text{for } x \in \mathbb{R}^N \text{ and }t>0, \\
  u(0,x)  =u_0(x) &\text{for } x \in \mathbb{R}^N,
    \end{array}
    \right.
\end{equation}
with data $u_0 \in L^1(\mathbb{R}^N)$ and exponents $0<s<1$ and $m>0$. Existence and uniqueness of a weak solution is established for
$m>m_c=(N-2s)_+/N$ giving rise to an $L^1$-contraction semigroup. Recently in \cite{DPQRV3}, it was proved the $C^{1,\alpha}$ regularity. Positivity of the
solution for any $m>0$ corresponding to non-negative data has been proved in \cite{BV2012}. We give a brief discussion on these facts in the Appendix.

%%%%%%%%%%%%%%%%%%%%%%%%%%%%%%%%%%%%%%%%%%%%
\subsection{Barenblatt Solutions of the Fractional Porous Medium Equation}\label{subsecBarenblatt}

An important tool that we use in the paper is represented by the so called Barenblatt solutions of the FPME. In \cite{VazquezBarenblattFractPME}, the second author proves existence, uniqueness and main properties of such fundamental solutions of the equation
\begin{equation}\label{eqPME}
u_t +(-\Delta)^{s}u^m=0,
\end{equation}
taking as initial data a Dirac delta $u(x,0)=M \delta(x),$ where $M>0$ is the mass of the solution. We will give here a short description of these functions and
recall their main properties we need in the paper. Next, we recall Theorem 1.1 from \cite{VazquezBarenblattFractPME}.

\begin{theorem}\label{ThUnicBarenblatt}
For every choice of parameters $s\in (0,1)$ and $m>m_c=\max\{(N-2s)/N,0\}$, and every $M>0$, Equation \eqref{eqPME} admits a unique fundamental solution; it is a
nonnegative and continuous weak solution for $t>0$ and takes the initial data in the sense of Radon measures. Such solution has the self-similar form
\begin{equation}\label{barenblattSol}
B_M(x,t)=t^{-\alpha}F_M(|x| t^{-\beta})
\end{equation}
for suitable $\alpha$ and $\beta$ that can be calculated in terms of $N$ and $s$ in a dimensional way, precisely
\begin{equation}\label{barenblattParam}
\alpha=\frac{N}{N(m-1)+2s}, \quad \beta=\frac{1}{N(m-1)+2s}.
\end{equation}
The profile function $F_M(r)$, $r>0$, is a bounded and H\"{o}lder continuous function, it is positive everywhere, it is monotone and goes to zero at infinity.
\end{theorem}

 By Theorem \ref{ThUnicBarenblatt} there exists a unique self-similar solution $B_1(x,t)$ with mass $M=1$ of Problem \eqref{eqPME}  and moreover, it has the form
$B_1(x,t)=t^{-\alpha} F_1(|x|t^{-\beta}).$ Let $B_{M}(x,t)$ the unique self-similar solution of Problem \eqref{eqPME} with mass $M$. Such function will be of the form
\begin{equation}\label{scalingBarenblatt}
B_{M}(x,t)=M B_1 \left(x, M^{m-1}t\right),
\end{equation}
which can be written in terms of the profile $F_1$ as
\begin{equation}\label{BM1}
B_{M}(x,t)=M^{1-(m-1)\alpha} t^{-\alpha} F\left( \left(M_1^{m-1}t \right)^{-\beta} |x|\right).
\end{equation}

\noindent Moreover, the precise characterization of the profile $F_M$ is given by Theorem $8.1$ of \cite{VazquezBarenblattFractPME}.

\begin{theorem}\label{ThProfileBarenblatt}
For every $m>m_1=N/(N+2s)$ we have the asymptotic estimate
\begin{equation}\label{decayF}
\lim_{r \rightarrow \infty}F_M(r) r^{N+2s}=C_1M^{\sigma},
\end{equation}
where $M=\int F(x)dx$, $C_1=C_1(m,N,s)>0$ and $\sigma=(m-m_1)(N+2s)\beta$. On the other hand, for $m_c<m<m_1$, there is a constant $C_{\infty}(m,N,s)$ such that
\begin{equation}\label{decayF2}
\lim_{r \rightarrow \infty}F_M(r) r^{2s/(1-m)}=C_{\infty}.
\end{equation}
\end{theorem}
The case $m=m_1$ has a logarithmic correction. The profile $F$ has the upper bound
\begin{equation}\label{decayF3}
F(r) \leq C r^{-N-2s+\epsilon}, \quad \forall r >0
\end{equation}
for every $\epsilon>0,$ and the lower bound
\begin{equation}\label{decayF4}
F(r) \geq C r^{-N-2s}\log r, \quad \text{ for all large }r.
\end{equation}

We state now some properties of the profile $F_1(r), \ r \geq 0$ obtained as consequences of formula \eqref{decayF} that we will use in what follows. Let us
consider first the case $m>m_1$.
\begin{enumerate}
  \item $F_1$ attains its maximum when $r=0$ i.e. $F(r) \leq F(0)$, for all $r \geq 0.$
  \item There exists $K_1>0$ such that
\begin{equation}\label{prop1F}
F(r) \leq K_1 r^{-(N+2s)}, \quad \forall r >0.
\end{equation}
  \item There exists $K_2>0$ such that
\begin{equation}\label{prop2F}
F(r) \geq K_2 (1+r^{N+2s})^{-1}, \quad \forall r \geq 0.
\end{equation}
\end{enumerate}
\noindent Similar estimates hold also in the case $m_c<m<m_1$, and the corresponding tail behavior is different, $F(r)\sim r^{-2s/(1-m)}.$ This will have an
effect in the different results we get for the generalized KPP problem.

As a consequence, the author also proves that the asymptotic behavior of general solutions of Problem \eqref{genfractPME} is represented by
such special solutions as described in Theorem 10.1 from \cite{VazquezBarenblattFractPME}.

\begin{theorem}\label{ThAsympGenfractPMEjlv}
 Let $u_0=\mu \in \mathcal{M}_+(\RN)$, let $M = \mu(\RN)$ and let $B_M$ be the self-similar Barenblatt solution with mass $M$. Then we have
 $$\lim_{t \rightarrow \infty}|u(x,t)-B_M(x,t;M)|=0$$
 and the convergence is uniform in $\RN$.
\end{theorem}

%%%%%%%%%%%%%%%%%%%%%%%%%%%%%%%%%%%%%%%%%%%%%%%%%%%%%%%%%%%%%%%%%%%%%%%%%%%%%%%%%%%%%%%%%%%%

\subsection{Lower estimates for nonnegative solutions in the case $m_1<m<1$}

\noindent We use the notations: $m_c:=(N-2s)/N$, $\vartheta:=1/[2s-N(1-m)]>0$ for $m>m_c$.
The results we quote are valid for initial data in a weighted space $u_0 \in L^1(\RR^N,\varphi dx)$, where $\varphi$ satisfies the following conditions:

\noindent \textbf{Assumption (A). } The function $\varphi \in C^2(\RN)$ is a positive real function that is radially symmetric and decreasing in $|x|\geq 1$.
Moreover $\varphi$ satisfies
$$0 \leq \varphi(x) \leq |x|^{-\alpha}\quad \text{for }|x|>>1 \text{ and }N-\frac{2s}{1-m} <\alpha <N+\frac{2s}{1-m}.$$

\noindent We recall now Theorem 4.1 from \cite{BV2012} giving local lower bounds for the solution of the diffusion problem.

\begin{theorem}[\textbf{Local lower bounds}]\label{Th41BV}
Let $R_0>0$, $m_c<m<1$ and let $0 \leq u_0 \in L^1(\RR^N,\varphi dx)$, where $\varphi$ is as in Assumption (A).  Let $u(\cdot,t) \in L^1(\RR^N,\varphi dx)$ be a
very weak solution to the Cauchy Problem \eqref{genfractPME},
corresponding to the initial datum $u_0$. Then there exists a time
\begin{equation}\label{tstar}
t_*:=C_* R_0^{\frac{1}{\vartheta}} \|u_0\|_{L^1(B_{R_0})}^{1-m}
\end{equation}
such that
\begin{equation}\label{LowerBoundSmallTimesFDE}
\inf_{x \in B_{R_0/2}}u(x,t) \geq K_1 R_0^{-\frac{2s}{1-m}}t^{\frac{1}{1-m}} \text{  if }0\leq t \leq t_*,
\end{equation}
and
\begin{equation}\label{LowerBoundLargeTimesFDE}
\inf_{x \in B_{R_0/2}}u(x,t) \geq K_1 \frac{\|u_0\|_{L^1(B_{R_0})}^{2s\vartheta}}{t^{N\vartheta}} \text{  if } t \geq t_*.
\end{equation}
The positive constants $C_*, \ K_1,\ K_2$ depend only on $m, \ s$ and $N\geq 1$.
\end{theorem}

\noindent The previous estimates, computed for $t=t_*$ rewrite as
\begin{equation}\label{LowerEstimtstar}
\inf_{x \in B_{R_0/2}}u(x,t) \geq K_1 C_*^{\frac{1}{1-m}}\|u_0\|_{L^1(B_{R_0})} R_0^{-N}.
\end{equation}
\noindent Then, if $R_0$ increases, the lower bound will decrease.

\noindent Concerning quantitative lower estimates, we recall Theorem 4.3 from \cite{BV2012}.

\begin{theorem}[\textbf{Global Lower Bounds when $m_1<m<1$}]\label{Th43BV}Under the conditions of Theorem \ref{Th41BV} we have in the range $m_1<m<1$
\begin{equation}\label{QuantLowerBoundFDE}
u(x,t) \geq \frac{C(t)}{|x|^{N+2s}} \text{   when } |x|>>1,
\end{equation}
valid for all $0<t<T$ with some bounded function $C>0$ that depends on $t,\ T$ and on the data.
\end{theorem}

\begin{theorem}[\textbf{Global Lower Bounds when $m_c<m<m_1$}]\label{ThGlobalLowerBounds1}
Under the conditions of Theorem \ref{Th41BV} we have in the range $m_c<m<m_1$
\begin{equation}
 u(x,t_0) \geq C(t) |x|^{-2s/(1-m)}
\end{equation}
if $|x|\geq R$ and $0<t<t_0.$

\end{theorem}

The lower estimates for exponents $m>1$ need a new analysis that we supply in the next section.

%%%%%%%%%%%%%%%%%%%%%%%%%%%%%%%%%%%%%%%%%%%%%%%%%%%%%%%%%%%%%%%%%%%%%%%%%%%%%%%%%%%%%%%%%%

\section{Lower parabolic estimate  in the case $m>1$}\label{SubsectLowerParEst}
We consider the FPME equation (with no reaction term) for $x\in\ren$ and $t>0$ with nonnegative and integrable initial data
\begin{equation}\label{eq.id}
  u(x,0)=u_0(x)\,,
\end{equation}
and we also assume that $u_0$ is bounded and has compact support or decays rapidly as $|x|\to\infty$.
We want to describe the  behavior of the solution \normalcolor $u(x,t)>0 $ as \normalcolor $|x|\to\infty$, more precisely its rate of decay, for small times $t>0$.  We take $m>1$ \normalcolor since the study of positivity for $m\le 1$ was dealt with in previous results.\normalcolor

The first step in our asymptotic positivity analysis of solutions of \eqref{eqPME}-\eqref{eq.id} is to ensure that solutions with positive data
remain positive in some region. We only need a special case that we quote next,
based on the positivity results of \cite{BV2012}.

\begin{theorem}[\textbf{Local lower bound}]\label{thm.lower.pme} Let $u$ be a weak solution to Equation \eqref{eqPME}, corresponding to $u_0\in L^1(\RR^d)$. Then
there exists a time
\begin{equation}\label{t*.PME}
t_*:=C \,R^{2s+N(m-1)}\,\|u_0\|_{L^1(B_{R})}^{-(m-1)}, \quad \vartheta:=1/[2s+N(m-1)]>0,
\end{equation}
such that for every $t\ge t_*$ we have the lower bound
\begin{equation}
\inf_{x\in B_{R/2}}u(x,t)\ge K\,\dfrac{\|u_0\|_{L^1(B_{R})}^{2s\vartheta}}{t^{N\vartheta}}
\end{equation}
 valid for all  $R>0$. The positive constants $C$ and $K$  depend only on $m,s$ and $N$, and
 not on $R$.
\end{theorem}

We may now state and prove the main lower estimate with precise tail behaviour, which is based on a delicate subsolution construction.

\begin{theorem}\label{low.par.est} Let $u(x,t)$ be a solution with initial data $u_0(x)\ge0$ such that $u_0(x)\ge 1$ in the ball $B_1(0)$. Then there is a time
$t_1>0$ and constants $C_*, R>0$ such that
\begin{equation}
u(x,t)\ge C_*\,t\,|x|^{-(N+2s)}
\end{equation}
if $|x|\ge R$  and $0<t<t_1$.
\end{theorem}

\noindent {\sl Proof.}  We consider the FPME for $m>1$ and initial data $u_0$  is  1 in the ball of radius 2. We use the previous result to prove that in the ball
of radius 1/2, then $u(x,t)\ge 1/2$ for $0<t<t_0$, a time that is calculated from the formulas above.

\medskip

\noindent $\bullet$ We want to construct a sub-solution of the form
$$
U^m(x,t)=G(|x|)+ t^m\,F^m(|x|).
$$
We want to choose $G\ge 0$ and $F\ge 0 $ in such a way that $U$ will be a formal sub-solution of the FPME in a domain of the form $Q=\{|x|\ge 1/2, 0<t<t_1\}$,
i.e., we want \ $U_t + (-\Delta)^s U^m\le0$ in $Q$.
Note  that
$$
U_t=(G(x)+ t^m\,F^m(x))^{(1/m)-1}t^{m-1}F^m(x)\le F(x).
$$
We also have, with $L_s=(-\Delta)^s$,
$$
L_sU^m=L_s G(x)+ t^m L_s F^m(x).
$$
We take $F$ positive, smooth and $F(r)\sim  r^{-(N+2s)}$ as $r\to\infty$ to get the desired conclusion after the comparison argument: $u(x,t)\ge U(x,t)\ge
ct\,r^{-(N+2s)}$ if  $r$ is large and  $t\sim 0$. For later use, let us say that  $F\le C_2r^{-(N+2s)}$ for $r>1/2$. Since $m>1$ we can choose $F$ smooth so  that
$ L_s F^m=O(r^{-(N+2s)})$ for $r>1/2$ (use the asymptotic estimates like the first lemma in \cite{BV2012})

We will take  $G(r)=0$ for $r=|x|\ge 1/2$ so that $U(x,t)=t\,F(x)$ there. If $G$ is also smooth we have $L_sG$ bounded and $L_sG\sim - C_1r^{-(N+2)}$ as
$r\to\infty$. By contracting $G$ in space, $\widetilde G(x)=G(kx)$, $k>0$, we may then say that $L_sG\le -C_1r^{-(N+2)}$ for $r>1/2$. Then we will have for
$r>1/2$
$$
U_t+ L_sU^m\le  F + L_sG+t^mL_s F^m\le  C_2r^{-(N+2s)} -C_1r^{-(N+2s)}+ t^m L_s F^m\le
$$
$$
(C_2+\varepsilon)r^{-(N+2s)} -C_1r^{-(N+2s)}\le 0
$$
if $C_1>C_2$. We can  choose $G$ large so that $C_1$ is large enough.

\medskip

\noindent $\bullet$ We now want to use  the viscosity method to compare $u(x,t)$ with $U(x,t)$ in the region $Q=\{|x|\ge 1/2, 0<t<t_1\}$, and this will prove that
$U(x,t)\le u(x,t)$ in $Q$.  Apart from the sub-solution condition that we have checked, we need suitable comparison of the boundary conditions at $r=1/2$
$$
U(1/2,t)=tF(1/2)\le 1/2.
$$
This ends the construction if the comparison result is justified. The contradiction argument at the first point of contact between $u$ and $U$ will be justified
as in \cite{BV2012} (where it was  applied to fast diffusion equations of fractional diffusion type) if the solution we have is a bit smooth: $u_t$ and $L_s u^m$
must be continuous and the equation must be satisfied pointwise there. This regularity is true and the proofs are under study now.

Alternatively, we may use Implicit Time Discretization with a sequence of approximations. The justification of the method in the elliptic case is done in the paper in collaboration with Volzone \cite{VazquezVolzone2012} on symmetrization techniques.  \finprf

\noindent {\bf Remark.} The level $u_0(x)\ge 1$ in the  ball $B_1(0)$ can be replaced by $u_0(x)\ge \varepsilon >0$ in any other ball by means
of translation and scaling. In this way the result is true for all continuous and nonnegative initial data $u_0$, of course nontrivial.

%%%%%%%%%%%%%%%%%%%%%%%%%%%%%%%%%%%%%%%%%%%%%%%%%%%%%%%%%%%%%%%%%%%%%%%%%%%%%%%%%%%%%%

\section{Evolution of level sets of solutions to Problem \eqref{KPP}}\label{SubsecConv0}

In this section we start the proof of the main result of the paper on evolution of level sets with exponential speed of propagation. In a first step we prove the
convergence to zero on outer sets. Since the decay assumption on the initial data is the same for $m_1<m<1$ and $m>1$, we will threat both cases, as well as
$m=1$, in the following lemma.

\begin{lemma}\label{lemmaConv0} We consider $m>m_1$ and let $u$ be the solution of Problem \eqref{KPP} with initial datum $u_0(x) \in L^1(\mathbb{R}^N)$, $0\leq
u_0\leq 1$. We assume that $u_0$ satisfies the decay property
\begin{equation}\label{dataDecay1}
u_0(x) \leq C|x|^{-(N+2s)}\quad  \text{for all } x\in \RN.
\end{equation}
Then, for $\sigma>\sigma_3$  if $m>1$ (respectively, for $\sigma> \sigma_{2}$ if $m_1<m\le1$), we have \
\begin{equation}\label{conv0}
u(x,t) \rightarrow 0 \quad \text{as}\quad t \rightarrow \infty
\end{equation}
uniformly for $|x|\geq e^{\sigma t}$.
\end{lemma}

\prf
We consider the solution $\overline{u}(x,t)$  of the linearized problem
$$
\overline{u}_t+L_s\overline{u}^m =f'(0) \overline{u}, \quad \overline{u}(0,x)=u_0(x).
$$
Since $f$ is a concave function, we have $f'(0)s\geq f(s)$ \ $ \forall s \in[0,1]$, and thus $\overline{u}$ is a super-solution of Problem \eqref{KPP}, which
implies the upper estimate
$$
u(x,t)\leq \overline{u}(x,t),\quad \text{  for  }t\geq 0, \ x \in \mathbb{R}^N.
$$
Next, we define $\tilde{v}(x,\tau)$ by
\begin{equation}\label{vtilde}
\tilde{v}(x,\tau)=e^{-f'(0) t}\overline{u}(x,t),
\end{equation}
and new time
%$$v_1(x)=\left\{
%  \begin{array}{ll}
%    a_1|x|^{-(N+2s)}, & \hbox{$|x|\geq \rho_1$;} \\
%    \epsilon=a_1 \rho_1^{-(N+2s)}, & \hbox{$|x|\leq \rho_1$.}
%  \end{array}
%\right.$$
\begin{equation}\label{tau1case}
\tau=\frac{1}{(m-1)f'(0)}\left[e^{(m-1)f'(0)t} -1 \right] \text{ if } m>1,\\
\end{equation}
\begin{equation}\label{tau2case}
\tau=\frac{1}{(1-m)f'(0)}\left[1-e^{-(1-m)f'(0)t}  \right] \text{ if } m<1,
\end{equation}
and $\tau=t$ for $m=1$. It is immediate to check that $\tilde{v}(x,\tau)$ is a solution of the FPME \eqref{FPME} with initial datum $\tilde{v}_0=u_0$. Let
$B_M(x,\tau)$ the Barenblatt solution with mass $M$ of the FPME, as defined in Section \ref{subsecBarenblatt}. By virtue of the properties of the Barenblatt
solutions and assumption \eqref{dataAssump} on the initial data, we conclude there exists $M>0$ big enough and $\tau_0>0$ such that
$$
\tilde{v}_0(x) \leq B_M(x,\tau_0).
$$
%otherwise we may start as initial data $\tilde{v}(\tau_1,\cdot)$, with a suitable choice of $\tau_1$.
By the Maximum Principle
$$
\tilde v(x,\tau) \leq B_M(x,\tau+\tau_0), \quad \forall x \in \mathbb{R}^N, t>0.
$$
Now, using the characterization of the decay of the Barenblatt profile given by relation \eqref{decayF}, we obtain that there exists $K_1>0$ such that
$F(r) \leq K_1r^{-(N+2s)},$ for all $r\geq0.$ We obtain the following upper estimate on the solution $u$ of Problem \eqref{KPP}:
\begin{align*}
u(x,t)&=e^{f'(0) t}\tilde{v}(x,\tau)\\
&\leq e^{f'(0) t} B_M(x,\tau+\tau_0)=e^{f'(0) t} (\tau+\tau_0)^{-\alpha}F(|x|(\tau+\tau_0)^{-\beta})\\
&\leq e^{f'(0) t} (\tau+\tau_0)^{-\alpha}K_1(|x|(\tau+\tau_0)^{-\beta}))^{-(N+2s)} \\
&=  K_1 e^{f'(0) t} (\tau+\tau_0)^{2\beta s}|x|^{-(N+2s)}.
\end{align*}
\noindent \textbf{Case $m>1$.} In order to continue the estimate, we remark that for large times $t$, the term $\tau^{2\beta s}$ has an influence on the result
only in the case $m>1$. Then $(\tau+\tau_0)^{2\beta s} \leq e^{(m-1)f'(0)t}$ for large $t$.
Let us assume that $|x|\geq e^{\sigma t}$. Then one has
\begin{align*}
u(x,t) \leq C K_1 e^{f'(0) t} \tau^{2\beta s}e^{-\sigma (N+2s)t} =C K_1 e^{[f'(0)  + 2f'(0)(m-1)\beta s  -\sigma (N+2s)]t}.
\end{align*}
We want to have  $f'(0)  + 2f'(0)(m-1)\beta s  -\sigma (N+2s)<0$, which is just the condition
$$
\sigma>\frac{1+2(m-1)\beta s}{N+2s}f'(0) =\sigma_3.
$$
We have obtained the convergence of $u(x,t)$ to $0$ as $t \rightarrow \infty$, for $|x|\geq e^{\sigma t}$.

\noindent \textbf{Case $m\le1$.} In this case, the term $(\tau+\tau_0)^{2\beta s}$ is bounded for every $t>0$ as we can see from \eqref{tau2case}. As before, we
assume $|x|\geq e^{\sigma t}$. Then, we get
\begin{align*}
u(x,t) \leq C K_1 e^{f'(0) t} e^{-\sigma (N+2s)t} =C K_1 e^{[f'(0)  -\sigma (N+2s)]t}.
\end{align*}
For $\sigma> \sigma_{2}=\frac{f'(0)}{N+2s}$, the exponent is negative $f'(0)  -\sigma (N+2s)<0$ and we obtain the convergence of $u(x,t)$ to $0$ as $t \rightarrow \infty$.

\finprf

\begin{lemma} We consider $m_c<m<m_1$. Let $u$ be the solution of problem \eqref{KPP} with initial datum $u_0(x) \in L^1(\mathbb{R}^N)$, $0\leq u_0\leq 1$ and we
assume $u_0$ satisfies the decay property
$$u_0(x) \leq C|x|^{-2s/(1-m)}\quad  \text{for all } x\in \RN.$$
Then, for $\sigma>\sigma_1$  we have
$$
u(x,t) \rightarrow 0, \quad t \rightarrow \infty
$$
uniformly for $|x|\geq e^{\sigma t}$.
\end{lemma}
\prf
The proof follows the same as in Lemma \ref{lemmaConv0} since the Barenblatt solution $B_M(x,\tau)=\tau^{-\alpha}F(|x|\tau^{-\beta})$ of the diffusion problem
satisfies $F(r) \sim r^{-2s/(1-m)}$ according to Theorem \ref{ThProfileBarenblatt}.
Therefore, we obtain the estimate
\begin{align*}
u(x,t)&\leq e^{f'(0) t} (\tau+\tau_0)^{-\alpha}K_1(|x|(\tau+\tau_0)^{-\beta}))^{-2s/(1-m)} \\
&=K_1 e^{f'(0) t}(\tau+\tau_0)^{-1/(1-m)}|x|^{-2s/(1-m)}.
\end{align*}
Since $m<1$, the term $(\tau+\tau_0)^{-1/(1-m)}$ is bounded and then, for $|x|\geq e^{\sigma t}$ we obtain
$$
u(x,t)\leq K_1 e^{f'(0) t -2s\sigma t /(1-m)}.
$$
For $\sigma> \frac{2s}{1-m}f'(0)=\sigma_1$ we obtain the desired convergence to $0$ as $t \rightarrow \infty.$

\finprf

\noindent\textbf{Remarks}

\noindent\textbf{I.} When $m=1$ we recover the minimal speed $\sigma_2=f'(0)/(N+2s)$ obtained by Cabr\'{e} and Roquejoffre in \cite{CabreRoquejoffreArxiv}. The
proof is similar, but in the nonlinear case we have to make an exponential change of time variable. Note also that we only use the decay properties of the
fundamental solution.

\noindent\textbf{II.} The value of the critical exponent $\sigma_2$ can be easily obtained from the following formal study of the level lines of $u(x,t)$.  Thus,
the set  $\{u(x,t) \sim \epsilon\}$
can be written in terms of $\tilde{v}(x,\tau)$ defined in \eqref{vtilde} as
\begin{equation}\label{vtilde.epsilon}
e^{f'(0) t}\tilde{v}(x,\tau) \sim \epsilon.
\end{equation}
By Theorem \ref{ThAsympGenfractPMEjlv}, $\tilde{v}(x,\tau)$ behaves like the Barenblatt solution of the Fractional Porous Medium Equation \eqref{genfractPME}:
$$
\tilde{v}(x,\tau) \sim B(x,\tau)=\tau^{-\alpha}F(r), \quad F(r)\sim r^{-(N+2s)}, \quad r=|x|\tau^{-\beta}.
$$
From \cite{VazquezBarenblattFractPME}, we know that
$
\displaystyle{B(x,\tau)\sim \tau^{-\alpha+\beta(N+2s)}|x|^{-(N+2s)}}
$, thus $\tilde{v}(x,\tau) \sim \tau^{2\beta s}|x|^{-(N+2s)}.$
At this moment, \eqref{vtilde.epsilon} implies
$\displaystyle{e^{f'(0) t} \tau^{2\beta s}|x|^{-(N+2s)}\sim \epsilon.}$

For instance in the $m>1$ case, it follows that
$$|x| \sim \left( \frac{1}{\epsilon} e^{ (1+2\beta s(m-1))f'(0) t} \right)^{1/(N+2s)}\sim  e^{ \frac{1+2\beta s(m-1)}{N+2s} f'(0) t} ,$$
and we deduce an exponential behavior of the level sets $|x| \sim  e^{\sigma_3 t},$ where $\displaystyle{\sigma_3=\frac{1+2\beta s(m-1)}{N+2s}f'(0).}$
Similarly, in the $m_1<m<1$ case, we get that
$$|x| \sim \left( \frac{1}{\epsilon} e^{ f'(0) t} \right)^{1/(N+2s)}\sim e^{\sigma_{2} t}, \quad  \sigma_{2}=\frac{f'(0)}{N+2s}.$$

%%%%%%%%%%%%%%%%%%%%%%%%%%%%%%%%%%%%%%%%%%%%%%%%%%%%%%%%%%%%%%%%%%%%%%%%%%%%%%

\section{Evolution of level sets II. Convergence to 1 on inner sets}\label{SubsecConv1}

In this section, we will prove the convergence to $1$ of the solution $u(x,t)$ of Problem \eqref{KPP}, i.\,e., the second part of the statements of our main theorems \ref{mainThm1}, \ref{mainThm3}, and \ref{mainThm3}.

\subsection{Case $m>m_1$ } We will present this case in full detail. The proof for the case $m_c<m<m_1$ being similar, we will sketch it at the end of this section. We have $N\geq 1$, $s \in (0,1)$, $m>m_1$, $f$ satisfies \eqref{propf}, and $\sigma_{2}=\frac{f'(0)}{N+2s}$ as defined in \eqref{sigma}.

\begin{proposition}\label{proConv1}
Let $N\geq 1$, $s\in (0,1)$, $m_1 < m$, $f$ satisfying \eqref{propf}. Let $u$ be a solution of Problem \eqref{KPP} with initial datum $0\leq u_0(\cdot)\leq 1$,
$u_0 \neq 0$. Then for every $\sigma \in (0,\sigma_{2})$, $u(x,t) \rightarrow 1$ uniformly on $\{|x|\leq e^{\sigma t}\}$ as $t \rightarrow \infty.$
\end{proposition}
\prf
We fix $\sigma \in (0,\sigma_{2})$. Proving the converge of $u(x,t)$ to $1$ is equivalent to proving the convergence of $1-u^m$ to $0$. Therefore, we fix $\lambda>0$ and we need to find a time
$t_\lambda$ large enough such that $1-u^m(x,t) \le \lambda$ for all $t\geq t_\lambda$ and $|x|\leq e^{\sigma t}.$

\noindent$\bullet$  Let us accept for the moment the following lower estimate that will be proved later as Lemma \ref{lemmaLowerBoundEps2}: given  $\nu \in (\sigma,  \sigma_{2})$, there exist $ \epsilon \in (0,1)$ and $t_0>0$ such that
\begin{equation}\label{uEps}
u\geq \epsilon, \quad \text{ for   } t \geq t_0 \text{   and   }|x|\leq  e^{\nu t}.
\end{equation}
 We now proceed with the last part of the argument, where the effect of the nonlinear diffusion is most clearly noticed. We take $t_1\geq t_0$ and consider the inner sets where
$$
\epsilon \le u \le 1  \quad \text{for } (x,t) \in \Omega_I:=\{ t \geq t_1, \ |x|\leq  e^{\nu t}\}.
$$
Let $v=1-u^m$. Then $\overline{v}$ satisfies the equation
\begin{equation}\label{eqvbar1}
\frac{1}{m}(1-v)^{\frac{1}{m} -1}v_t  +L_s v+ f(u)=0,
\end{equation}
that we write in the form
\begin{equation}
\label{eqvbar}a(x,t)v_t  +L_s \overline{v}+ b(x,t)v=0, \qquad a(x,t)=\frac{1}{m} u^{1-m}, \quad b(x,t)=\frac{f(u)}{v}.
\end{equation}
Moreover, we estimate $a(x,t)$ as follows
$$
  a_0 = \frac{1}{m} \epsilon^{1-m} \leq a(x,t) \leq a_1:=\frac{1}{m} \text{ in }\Omega_I, \ \text{ if }\ m<1,
$$
respectively,
$$
  a_0 = \frac{1}{m}  \leq a(x,t) \leq a_1:=\frac{1}{m}\epsilon^{1-m} \text{ in }\Omega_I, \ \text{ if } m>1.
$$
We argue similarly for $b(x,t)$ in $\Omega_I$:
$$
b(x,t)=\frac{f(u)}{1-u^m} = \frac{f(u)}{(1-u) m \xi^{m-1}}\geq b_0, \quad \xi \in (u,1),
$$
where
$$
b_0=\frac{1}{m}\frac{f(\epsilon)}{1-\epsilon} \epsilon^{1-m} \text{ if } m<1 \quad \text{and}\quad b_0=\frac{1}{m}\frac{f(\epsilon)}{1-\epsilon}  \text{ if }m>1.
$$
In particular, $\overline{v}$ satisfies
\begin{equation}\label{eqv2}
a(x,t)v_t  +L_s v+ b_0v\leq0 \quad \text{in }\Omega_I.
\end{equation}

\medskip

\noindent$\bullet$ We look for  a super-solution $w$ to Problem \eqref{eqvbar} that will be found as a solution to a linear problem with constant coefficients, and we also need  that $w_t\le 0$. More precisely, we consider $w$ solution of the concrete problem
\begin{equation}\label{eqw2}
  \left\{ \begin{array}{ll}
  a_1 w_{t}(x,t) + L_s w(x,t)+b_0 w=0  &\text{for } x \in \mathbb{R}^N \text{ and }t>0, \\[2mm]
 w(x,t_1)  =1+\frac{1}{C_2}|x|^{\gamma}&\text{for } x \in \mathbb{R}^N.
    \end{array}
    \right.
\end{equation}
where the exponent $\gamma$ taken such that
\begin{equation}\label{gamma}
0<\gamma:=\frac{1}{\nu}\frac{b_0}{a_1} <2s.
\end{equation}
We can eventual consider a smaller $\epsilon$ for this inequality to hold. Equation \eqref{eqw2} is linear, the solution can be computed explicitly
$$
w(x,t)=e^{-\frac{b_0}{a_1}(t-t_1)} \overline{w} (x,\tau), \qquad \tau=\frac{1}{a_1}(t-t_1),
$$
where $\overline{w}(x,\tau)$ solves the linear problem
$$
\overline{w}_{\tau}(x,\tau) +L_s \overline{w}(x,\tau)=0, \qquad \overline{w}(0)=1+\frac{1}{C_2 }|x|^{\gamma}.
$$
We observe that $\tilde{w}$ can be written in the following form
\begin{equation}\label{wbar}
\overline{w}(x,\tau)= 1+\frac{1}{ C_2} U(x,\tau+\tau_1),
\end{equation}where
$$U(x,\tau)=\tau^{\alpha_1}F(|x|\tau^{-\beta_1}), \quad \alpha_1=\frac{\gamma}{2s}, \ \beta_1=\frac{1}{2s},$$
is the self-similar solution of the linear problem
$$
U_{\tau} (x,\tau) + L_s U(x,\tau)=0, \qquad U(x,0)=|x|^{\gamma}.
$$
The properties of the self-similar solutions $U(x,\tau)$ deserve a separate study, which is done in detail in Section \ref{sectSelfSimLinProb}. Thus, by Lemma \ref{lemmaSelfSim} the profile $F$ is non-decreasing and $U(x,\tau)$ has a spatial decay as $|x|^{\gamma}$ for large $|x|\tau^{-1/2s}$:
%$$
%\eta F'(\eta) \ge C_3 F(\eta), \quad \forall \eta >0.
%$$
%$$F(\eta) \sim \eta^{\gamma} \quad \text{for large  }\eta,$$
%which implies
\begin{equation}\label{U}
C_2 |x|^{\gamma}\le U(x,\tau) \le C_1 |x|^{\gamma} \quad \text{for all} \quad |x|\tau^{-1/2s}\geq K_1.
\end{equation}
We will consider a suitable delay time $\tau_1$ in the definition of $\overline{w}$ stated in \eqref{wbar}.
In what follows we will use the notation $\eta=|x|\tau^{-\beta_1}$. We check that the derivative $w_t$ is negative:
\begin{align*}w_t(x,t)&=\frac{d}{dt} \left[ e^{-\frac{b_0}{a_1}(t-t_1)} (1+ C_2^{-1}U(x,\tau+\tau_1) ) \right] \\
&=e^{-\frac{b_0}{a_1}(t-t_1)}  \left[
-\frac{b_0}{a_1}(1+ C_2^{-1}(\tau+\tau_1)^{\alpha_1}F(\eta) ) + C_2^{-1}(\tau+\tau_1)^{\alpha_1-1}\left(\alpha_1  F(\eta)- \beta_1 \eta F'(\eta) \right)
 \frac{d\tau}{dt}\right]\\
&=e^{-\frac{b_0}{a_1}(t-t_1)} \frac{1}{ a_1C_2 } \left[
-b_0C_2  + ( -b_0 (\tau+\tau_1) + \alpha_1 )(\tau+\tau_1)^{\alpha_1-1}F(\eta)  - \beta_1(\tau+\tau_1)^{\alpha_1-1}  \eta F'(\eta)
 \right]
\end{align*}
Since $F'(\eta)>0$ for all $\eta>0$, we get that $w_t(x,t) \leq 0$ for all $t\ge t_1$ if $\tau + \tau_1 \ge \alpha_1/b_0,$ which is true for a suitable choice of $\tau_1$.
%In terms of $t$, this means that $\displaystyle{t \geq t_1+ a_1 \frac{\alpha_1}{b_0}= t_1+ \frac{1}{2s\nu}}=:t_2.$

%\noindent \textbf{Remark.} We continue our analysis for times $t>t_2$. Therefore, we consider the set $\Omega'_I=\{ t \geq t_2, \ |x|\leq  e^{\nu t}\}$ which is a
%subset of $\Omega_I$ previously defined.
\medskip

\noindent$\bullet$ Now we can compare $w$ and $v$ by applying the Maximum Principle stated in Lemma \ref{MaxPrinc} of the Appendix , as in \cite{CabreRoquejoffreArxiv}. Define $\overline{v}=v-w$ and ensure the hypothesis of the Lemma are satisfied.

(H1) We check that $w(x,t_1)\ge v(x,t_1)$ for all $x \in \RN$:
\begin{align*}
w(x,t_1) &\ge 1 > v=1-u^m, \quad \forall x \in \RN.
\end{align*}

(H2) We check that $w\ge v$ in $\left((t_1,\infty)\times \RN \right)\setminus  \Omega_I$, that is $t\geq t_1$ and $|x|\ge e^{\nu t}.$ At this point, we use the estimates \eqref{U}. We ensure that $e^{\nu t}\ge K_1 (\tau+\tau_1)^{1/2s}$ for all $t\geq t_1$, which is true by choosing eventually a larger $t_1$. Therefore
\begin{align*}
w(x,t) &=e^{-\frac{b_0}{a_1}(t-t_1)} \tilde{w} \left(\tau,x\right) \ge e^{-\frac{b_0}{a_1}(t-t_1)} (1+  \frac{1}{C_2 }
C_2 |x|^{\gamma} ) \\
&\ge
e^{-\frac{b_0}{a_1}(t-t_1)}(1+ e^{\gamma \nu t})
 \geq 1 \ge v(x,t) \quad \text{for all  }t\geq t_1, \ |x|\ge e^{\nu t}
\end{align*}
since $\gamma$ satisfies \eqref{gamma}. By the previous computation $\overline{v} \leq 0$ in $\left((t_1,\infty)\times \RN
\right)\setminus  \Omega_I$.

(H3) Next step is to prove that $\overline{v}$ is a sub-solution of Problem
\eqref{eqv2}. Indeed, we have that
$$
a(x,t)\overline{v}_t + L_s \overline{v} +b_0 \overline{v}= a(x,t)v_t+ L_sv+b_0v - \left[a_1w_t +L_s w+b_0w \right] + (a_1-a(x,t))w_t \leq 0 \text{  in
}\Omega_I.
$$
By Lemma \ref{MaxPrinc} we obtain that $\overline{v} \leq 0 $ in $[t_1,\infty)\times \RN$ for $t_1$ taken to be
large enough. Thus,
\begin{align*}
v(x,t) \le w(x,t)= e^{-{\frac{b_0}{a_1}}(t-t_1)} (1+ C_2^{-1}U(x,\tau))   \le  e^{-{\frac{b_0}{a_1}}(t-t_1)} (1+ \frac{C_1}{C_2}  |x|^{\gamma}).
\end{align*}
%We consider the constants $C_\lambda$ small enough such that
%$$
%e^{\frac{b_0}{a_1}t_1}\frac{C_1}{C_2} C_\lambda^\gamma <\lambda.
%$$
%Once we choose $C_\lambda$, we take $ t_\lambda$ large enough such that
%$$
%1 \leq C_\lambda e^{(\sigma -\nu)t_{\lambda}} \quad  \text{and}\quad e^{-{\frac{b_0}{a_1}}(t_\lambda-t_1)} \le \frac{\lambda}{2}.
%$$
%Observe that the time $t_\lambda$ depends only on $\sigma,\nu, a_0,b_0,t_1$ and not on the point $|x|$.

\noindent$\bullet$ Let us consider the inner set $(x,t) \in \{ t\geq t_{\lambda},\ |x|\leq C_{\lambda} e^{\nu t}\}$. We have
\begin{align*}
v(x,t) &\leq   e^{-{\frac{b_0}{a_1}}(t-t_1)} (1+ \frac{C_1}{C_2}C_{\lambda}^{\gamma} e^{\gamma \nu t})\le
 e^{-{\frac{b_0}{a_1}}(t_\lambda-t_1)} + \frac{C_1}{C_2} e^{\frac{b_0}{a_1}t_1} C_{\lambda}^{\gamma}  \le \lambda
\end{align*}
for $C_\lambda$ small enough and $t_\lambda$ large enough.

Finally, since $\sigma<\nu $ then $e^{\sigma t} \le C_{\lambda} e^{\nu t} $ for every $t\geq t_{\lambda}$ with $t_{\lambda}$ large enough, and the previous
inequality implies that
$$1-u^m(x,t)=v(x,t) \leq \lambda  \qquad \text{for}\quad t\geq t_{\lambda}, \ |x|\leq e^{\sigma t}, $$
which concluded the proof of the uniform convergence to the level $u=1$. \finprf

%%%%%%%%%%%%%%%%%%%%%%%%%%%%%%%%%%%%%%%%%%%%%%%%%%%%%

To complete the proof of the result of this subsection, we need to supply the proof of the lower estimate \eqref{uEps}. This will be done in three steps,

\smallskip

\noindent Step I. Starting with arbitrary initial datum $0\leq u_0\le 1$, $u_0 \neq 0,$  we obtain a lower bound for $u$ with the desired tail $u\ge c\, |x|^{-(N+2s)}$ for large $|x|$. The result corresponds to Lemma \ref{lemmaTail}.

\noindent Step II. We prove that given an initial data taking the value $\epsilon$ in the ball of
radius $\rho_0$ and decaying like that $|x|^{-(N+2s)}$ for large $|x|$, the corresponding solution of Problem \eqref{KPP} will be raised to at least the same level $\epsilon$ in a larger ball $\rho_1$ and in a later time that is estimated. The sizes are important. This will be Lemma \ref{Lemmav0}.

\noindent Step III. By combining the previous two results, we conclude that $u\geq \epsilon$ on the inner sets, for a certain $\epsilon>0$.
This will be Lemma \ref{lemmaLowerBoundEps} and Lemma \ref{lemmaLowerBoundEps2}.

Steps II and III follow the ideas of \cite{CabreRoquejoffreArxiv} in the linear case, with a long technical adaptation to nonlinear diffusion.

\smallskip

\begin{lemma}[\bf Long Tail Behaviour]\label{lemmaTail}
Let $N\geq 1$, $s\in (0,1)$, $m>m_1$, $f$ satisfying \eqref{propf} and $\sigma \in (0,\sigma_{2})$. Let $u$ be the solution of Problem \eqref{KPP} with initial datum $u(0,\cdot)=u_0$, where $0\leq u_0 \leq 1$, $u_0 \neq 0$. Then for any fixed $t_0>0$ there exist $\epsilon \in (0,1)$, $a_0>0$, $\rho_0>1$ such that
$$u(x,t) \geq v_0(x):=\left\{
  \begin{array}{ll}
    a_0|x|^{-(N+2s)}, & \hbox{$|x|\geq \rho_0$,} \\
    \epsilon=a_0 \rho_0^{-(N+2s)}, & \hbox{$|x|\leq \rho_0$,}
  \end{array}
\right.$$
for all $t \in [t_0,2t_0]$, $x\in \RN.$
\end{lemma}
\prf We recall that $\sigma_2=f'(0)/(N+2s)$. The idea is to view $u$ the solution of Problem \eqref{KPP} as a super-solution of the homogeneous problem with the same initial datum $u_0$, that is the FPME. Therefore,
$$u(x,t_0 +t) \geq \underline{u}(x,t ), \quad \forall t \geq 0, \quad x \in \RN,$$
where $\underline{u}$ is the solution of the FPME with initial datum $u_0$
\begin{equation}\label{fracFDE}
  \left\{ \begin{array}{ll}
  \underline{u}_{t}(x,t) + L_s \underline{u}^m(x,t)=0 &\text{for } x \in \mathbb{R}^N \text{ and }t>0, \\
 \underline{u}(x,0)  =u_0(x) &\text{for } x \in \mathbb{R}^N.
    \end{array}
    \right.
\end{equation}
We will estimate $\underline{u}$ from below by using the local and global estimates on the FPME given in Theorem \ref{Th41BV} and Theorem \ref{Th43BV} for $m<1$,
respectively Theorem \ref{thm.lower.pme} for $m>1$. The decay in case $m=1$ is well known, see Section \ref{sectSelfSimLinProb} for a review.
More exactly, in all cases $m>m_1$, there exist a time $T>0$ and constant $R>0$ such that
$$
u(x,t) \geq C(t) |x|^{-(N+2s)}, \quad \forall |x|\geq R, \ 0<t< T.
$$
Then, for a fixed $t_* \in (0,T)$ which also satisfies $t_*<t_0$, we can find a Barenblatt solution $B_{M}(x,t)$ and a time $t_2>0$ such that
$$
u(x,t_*) \geq B_M(x,t_2), \quad \forall x \in \RN,
$$
and therefore, by the Comparison Principle
$$
u(x,t+t_*) \geq B_M(x,t+t_2), \quad \forall x \in \RN, \ t \geq 0.
$$
In particular, we can choose $\epsilon>0$ such that
$$
u(x,t) \geq v_0(x):=\left\{
  \begin{array}{ll}
    a_0|x|^{-(N+2s)}, & \hbox{$|x|\geq \rho_0$,} \\
    \epsilon=a_0 \rho_0^{-(N+2s)}, & \hbox{$|x|\leq \rho_0$,}
  \end{array}
\right.$$
for all $x \in \RN, \ t \in [t_0,2t_0].$

\finprf

\begin{lemma}[\bf Positivity for a sequence of times]\label{Lemmav0} Let $m>m_1$. For every $\sigma < \sigma_{2}$ there exist $t_0 \geq 1$ and $0<\epsilon_0<1$ depending only on $N,s,f$ and $\sigma$ for which the following holds: given
$\rho_0\geq 1$ and $0<\epsilon \leq \epsilon_0$, let $a_0>0$ be defined by $a_0 \rho_0^{-(N+2s)}=\epsilon$, if we take
\begin{equation}\label{v0def}
v_0(x)=\left\{
  \begin{array}{ll}
    a_0|x|^{-(N+2s)}, & \hbox{$|x|\geq \rho_0$,} \\
    \epsilon=a_0 \rho_0^{-(N+2s)}, & \hbox{$|x|\leq \rho_0$\,,}
  \end{array}
\right.\end{equation}
then the solution $v$ of Problem \eqref{KPP} with initial condition $v_0$ satisfies
\begin{equation}\label{vxkt0}
v(x,kt_0) \geq \epsilon \quad \text{for }|x|\leq \rho_0e^{\sigma k t_0},
\end{equation}
for all $k\in \{0,1,2,3,...\}.$
\end{lemma}

\prf \noindent {\sc Case $m\ge1$.}

I. {\sl Preliminary choices}. From the beginning we fix $\sigma \in (0,\sigma_{2})$. We will do a very detailed analysis of the case $k=1$, which is then iterated for the rest of values of $k$.
We take $\delta \in (0,1)$ small enough such that
\begin{equation}\label{delta}
\frac{f(\delta)}{(N+2s)\delta} \geq \sigma, \quad \frac{f(\delta)}{(N+2s)\delta} \geq N(m-1)\beta \sigma_2.
\end{equation}
For example, take $\delta$ such that
$$
\frac{f(\delta)}{(N+2s)\delta}=\frac{1}{2}\left( \sigma_2 + \max\{\sigma, N(m-1)\beta \sigma_2 \}\right).
$$
This choice will be explained later.  Next we take $t_0$ sufficiently large depending only on $N,s,u_0$ and $\sigma$ such that
\begin{equation}\label{condt0}
  e^{\frac{f(\delta)}{\delta}t_0}\left(1+( \tau_0/\tau_1(1))\right)^{-N \beta}\ge K_3,  \qquad (K_2/2K_1)^{1/(N+2s)}\, e^{\frac{f(\delta)}{(N+2s)\delta}t_0} \ge e^{\sigma t_0}.
\end{equation}
where $\tau_1(1)=\epsilon^{1-m}c_2$ with $c_2$ a positive constant that we state explicitly later, and $K_2<2K_1$ are constants describing the properties of the profile $F_1$ of the Barenblatt function with mass $1$ given in \eqref{prop1F} and \eqref{prop2F}, and we recall for convenience that
$$
K_2 (1+r^{N+2s})^{-1} \le F(r) \le K_1 r^{-(N+2s)}, \qquad \forall r >0.
$$
Define now $\epsilon_0 \in (0,\sigma)$ by  $$
\epsilon_0=\delta e^{-f'(0)t_0} .
$$
Clearly, $\epsilon_0<\delta$. Now, we fix $0<\epsilon<\epsilon_0$ and $\rho_0 >1$.

\smallskip

\noindent II. {\sl Construction of sub-solutions to Problem \eqref{KPP}}. Let $w$ be a solution of the  problem with linearized reaction
\begin{equation}\label{eqw}
  \left\{ \begin{array}{ll}
  w_{t}(x,t) + L_s w^m(x,t)=\frac{f(\delta)}{\delta}w &\text{for } x \in \mathbb{R}^N \text{ and }t>0, \\
  w(0,x)  =v_0(x) &\text{for } x \in \mathbb{R}^N.
    \end{array}
    \right.
\end{equation}
We define $\overline{w}(x,\tau)$ by
$$\overline{w}(x,\tau)=e^{-\frac{f(\delta)}{\delta}t} w(x,t),$$
with a new time
\begin{equation}\label{tauDef1}
\tau=\frac{1}{(m-1)f'(0)}\left[e^{(m-1)f'(0)t} -1 \right] \text{ if } m>1,
\end{equation}
so that $\tau=t$ in the limit $m=1$. Then, $\overline{w}$ is a solution of the Fractional Porous Medium Equation with initial datum $v_0$
\begin{equation}\label{PMEv0}
  \left\{ \begin{array}{ll}
  \overline{w}_{\tau}(x,\tau) + L_s \overline{w}^m(x,\tau)=0 &\text{for } x \in \mathbb{R}^N \text{ and }\tau>0, \\[1mm]
  \overline{w}(x,0)  =v_0(x) &\text{for } x \in \mathbb{R}^N.
    \end{array}
    \right.
\end{equation}

%%%%%%%%%%%%%%%%%%%%%%%%%%%%%%%%%%%%%%%%%%%%%%%%%%%%%%%%%%%%%%%%%\hrule

\medskip

\noindent III. {\sl Comparison with a Barenblatt solution. Lower bound for $v(x,t_0)$.} We prove that there exist $M_1>0$ and $\tau_1>0$ such that
\begin{equation}\label{v0BM1}
v_0(x) \geq B_{M_1}(x,\tau_1), \quad \forall x \in \RN,
\end{equation}
where $B_{M_1}(x,\tau)$ is the Barenblatt solution of Problem \eqref{FPME} with mass $M_1$ given by Theorem \ref{ThUnicBarenblatt}:
\begin{equation}\label{scalingBarenblatt}
B_{M_1}(x,\tau)=M_1 B_1 \left( x, M_1^{m-1}\tau\right).
\end{equation}
Now, $B_{M_1}(x,\tau)$ can be written in terms of the profile $F_1$ as
\begin{equation}\label{BM1}
B_{M_1}(x,\tau)=M_1^{1-(m-1)\alpha} \tau^{-\alpha} F_1\left( \left(M_1^{m-1}\tau \right)^{-\beta} |x|\right).
\end{equation}
We will use the properties of the profile $F_1$ stated in \eqref{prop1F} and \eqref{prop2F}. With this information, we will find the constants $M_1>0$ and $\tau_1>0$ such that inequality \eqref{v0BM1} at the initial time holds true. For $|x|\leq \rho_0$ we have that
$B_{M_1}(x,\tau_1)\leq M_1^{1-(m-1)\alpha} \tau_1^{-\alpha} F_1(0).$ Note that $1-(m-1)\alpha=2\beta s>0$. We impose the first condition
\begin{equation}\label{cond1}
M_1^{2\beta s} \tau_1^{-\alpha} F_1(0)\leq \epsilon.
\end{equation}
Let $|x|\geq \rho_0$. Then,
$$
B_{M_1}(x,\tau)\leq M_1^{2\beta s} \tau_1^{-\alpha}K_1  \left( \left(M_1^{m-1}\tau \right)^{-\beta} |x|\right)^{-(N+2s)}.
$$
In order to use this the inequality for large $|x|$ we also impose the condition
\begin{equation}\label{cond2}
K_1M_1^{1+2\beta s (m-1)}\tau_1^{2 \beta s}\le a_0, \quad \text{where } a_0= \epsilon \rho_0^{N+2s}.
\end{equation}
Conditions \eqref{cond1} and \eqref{cond2} are sufficient for inequality \eqref{v0BM1} to hold.
Under such restrictions the stated inequality \eqref{v0BM1} holds true. Then, by the Comparison Principle we get
\begin{equation}\label{BM1below}
B_{M_1}(x,\tau+\tau_1) \leq \overline{w}(x,\tau),\quad \text{for all } |x| \in \RN, \ \tau>0.
\end{equation}
Putting equality in the inequalities  \eqref{cond1} and \eqref{cond2}  we get
\begin{equation}\label{cond3}
M_1=c_1 \ve\rho_0^{N}, \qquad \tau_1=c_2\ve^{1-m}\rho_0^{2s}\,,
\end{equation}
(with $c_1, c_2$ positive constants not depending on $\ve$ or $\rho$). We can easily see that the expressions are dimensionally correct.
The constants are given by  $c_1=F_1(0)^{-2s/(N+2s)}K_1^{-N/(N+2s)},$ $c_2=F_1(0)^{(1+2(m-1)\beta s)/(\beta(N+2s))}K_1^{-2s/(N+2s)}$. In particular, $(M^{m-1}\tau_1)^{\beta}=c_3\rho_0$, with $c_3=c_1^{(m-1)\beta}c_2^{\beta}.$

Since $v_0\leq \epsilon$ in $\RN$ then $\overline{w}(x,\tau) \leq \epsilon$ for all $x \in \RN$, $\tau >0$, and then in terms of $w(x,t)$ we obtain the following
bound
$$
0 \leq w(x,t) \leq e^{\frac{f(\delta)}{\delta}t_0} \epsilon \leq e^{f'(0)t_0} \epsilon_0 =\delta, \quad \forall t \leq t_0.
$$
Since $\frac{f(\delta)}{\delta} w \leq f(w)$ for $0\leq w \leq \delta$,  $w$ is a sub-solution of Problem \eqref{KPP} in $[0,t_0] \times \RN$. By the Comparison Principle and estimate \eqref{BM1below} we obtain that at the moment $t_0$
\begin{equation}\label{lowert0}
v(\cdot,t_0) \geq w(\cdot,t_0) = e^{\frac{f(\delta)}{\delta}t_0}\overline{w}(\cdot,\tau_0 ) \geq
e^{\frac{f(\delta)}{\delta}t_0}B_{M_1}(\cdot,\tau_0+\tau_1 ) \ \text{ in } \RN,
\end{equation}
where we use the notation $\tau_0=\tau(t_0)$ defined by \eqref{tauDef1}.

\medskip

\noindent IV. We will now prove that estimate \eqref{lowert0} with the choices \eqref{cond3} for $M_1$ and $\tau_1$ implies the lower bound stated in Lemma \ref{Lemmav0} in the case $k=1$, $m>1$.
Indeed,  we have (*)
\begin{align*}
v(x, t_0) &\geq  e^{\frac{f(\delta)}{\delta}t_0}B_{M_1}(x,\tau_0+\tau_1 )  \\
&=e^{\frac{f(\delta)}{\delta}t_0}M_1^{2\beta s}( \tau_0+\tau_1)^{-\alpha}F_1 \left( M_1^{-(m-1)\beta}(\tau_0+\tau_1)^{-\beta} |x|  \right) \\
&\geq e^{\frac{f(\delta)}{\delta}t_0} M_1^{2\beta s}( \tau_0+\tau_1)^{-\alpha}K_2  \left( 1+ (M_1^{-(m-1)\beta}(\tau_0+\tau_1)^{-\beta}
|x|)^{(N+2s)}\right)^{-1}\\
&= K_2 e^{\frac{f(\delta)}{\delta}t_0} M_1^{2\beta s}( \tau_0+\tau_1)^{-\alpha}  \left( 1+ c_3^{-(N+2s)}\left(1+( \tau_0/\tau_1)\right)^{-\beta (N+2s)}
(|x|/\rho_0)^{(N+2s)}\right)^{-1}
\end{align*}
Our aim now is to be able to continue this estimate until we reach a bound $v_1(x)$ of the form \eqref{v0def} for the same $\epsilon$ and a different  parameter $\rho_1$.
We will choose some $\rho_1$ and then check that the lower bound for $v(x,t_0)$ is larger than $\epsilon$ at $|x|=\rho_1$. In order to simplify the estimate of the last parenthesis in formula (*), we will impose the condition
$$
c_3^{-(N+2s)} \left(1+( \tau_0/\tau_1)\right)^{-\beta (N+2s)}
(\rho_1/\rho_0)^{(N+2s)}\ge 1\,,
$$
and then we only need to have
\begin{equation}\label{vrho1t0}
v(\rho_1,t_0)\ge (K_2/2) e^{\frac{f(\delta)}{\delta}t_0} M_1^{1+2(m-1)\beta s}( \tau_0+\tau_1)^{2\beta s} \rho_1^{-(N+2s)}\ge \epsilon.
\end{equation}
The first condition is equivalent to
$$
(\rho_1/\rho_0)^{(N+2s)}\ge c_3^{N+2s}\left(1+( \tau_0/\tau_1)\right)^{\beta (N+2s)}
$$
while, taking into account that $M_1^{1+2(m-1)\beta s}\tau_1^{2\beta s}=a_0/K_1$ and $a_0=\epsilon \rho_0^{N+2s}$, the second means that
\begin{equation}\label{rho1.upper}
(\rho_1/\rho_0)^{(N+2s)}\le (K_2/2K_1) e^{\frac{f(\delta)}{\delta}t_0}\left(1+( \tau_0/\tau_1)\right)^{2s\beta}\,.
\end{equation}
Both conditions are compatible iff
\begin{equation*}\label{t0def1}
e^{\frac{f(\delta)}{\delta}t_0}\left(1+( \tau_0/\tau_1)\right)^{-N \beta}\ge \frac{2K_1}{K_2}c_3^{N+2s}=K_3.
\end{equation*}
Now recall that $\tau_1$ depends on $\rho_0$ by \eqref{cond3}, and $\tau_1$ is bounded below by $\tau_1(1)$, the value for $\rho=1$. We see this condition as a way of choosing $t_0$. Using the fact that for large $t$ we have $$\tau\sim e^{(m-1)f'(0)t}/(m-1)f'(0),
$$
we easily see that for large $t_0$ the left-hand side looks like
$$
C\,e^{\frac{f(\delta)}{\delta}t_0-N(m-1)\beta f'(0)t_0},
$$
hence, the compatibility condition can be solved if \ $ f(\delta)/\delta> N(m-1)\beta f'(0)$. Since $\delta$ is small enough so that $ f(\delta)/\delta\approx f'(0)$, this means that we need $N(m-1)\beta <1$ which is true. We conclude that there exists $t_0>0$ large enough such that
\begin{equation}\label{t0def1}
e^{\frac{f(\delta)}{\delta}t_0}\left(1+( \tau_0/\tau_1(1))\right)^{-N \beta}\ge K_3\,.
\end{equation}
 This choice of $t_0$ is independent of $\rho_0$.

 Once this is guaranteed, we choose the largest possible $\rho_1$ satisfying \eqref{rho1.upper}, which is
\begin{equation}\label{rho1}
\frac{\rho_1}{\rho_0}= (K_2/2K_1)^{1/(N+2s)}\, e^{\frac{f(\delta)}{(N+2s)\delta}t_0}\left(1+( \tau_0/\tau_1)\right)^{2s\beta/(N+2s)}:=L_0.
\end{equation}

\medskip

\noindent V. With this choice of $\rho_1$ and $t_0$,  estimate \eqref{vrho1t0} holds. Going back to Point IV above, we have
$$
v(x,t_0)\geq e^{\frac{f(\delta)}{\delta}t_0}B_{M_1}(x,\tau(t_0)+\tau_1 )\geq \epsilon, \quad \text{ for }|x|=\rho_1\,,
$$
and thus, since the profile $F_1$ is non-increasing we get that
$$
v(x,t_0)\geq \epsilon,\quad \forall |x|\leq \rho_1\,.
$$
Using \eqref{rho1} we get that
$$
v(x,t_0)\geq \epsilon \rho_1^{N+2s} |x|^{-(N+2s)}, \quad \forall |x|\geq \rho_1.
$$
Finally, we define $a_1:=\epsilon \rho_1^{N+2s}$ and thus $v(\cdot,t_0) \geq v_1(\cdot)$ where $v_1$ is given by the expression
$$v_1(x)=\left\{
  \begin{array}{ll}
    a_1|x|^{-(N+2s)}, & \hbox{$|x|\geq \rho_1$;} \\
    \epsilon=a_1 \rho_1^{-(N+2s)}, & \hbox{$|x|\leq \rho_1$.}
  \end{array}
\right.$$

\begin{figure}[h!]
  \centering
  \includegraphics[width=120mm,height=45mm]{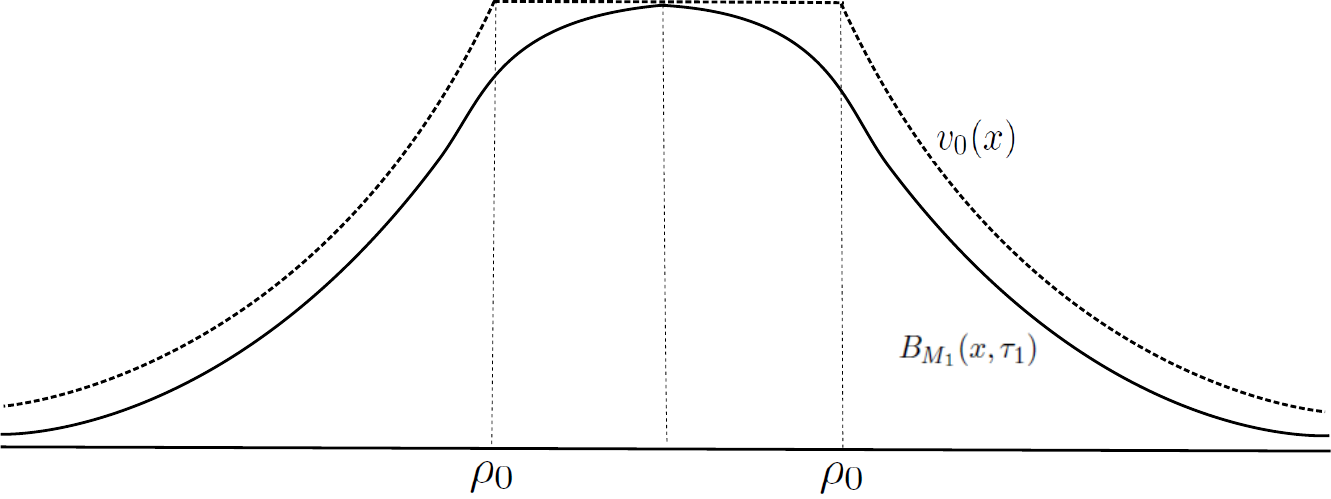}\caption{Step $III$}\label{figv0}
  \includegraphics[width=120mm,height=45mm]{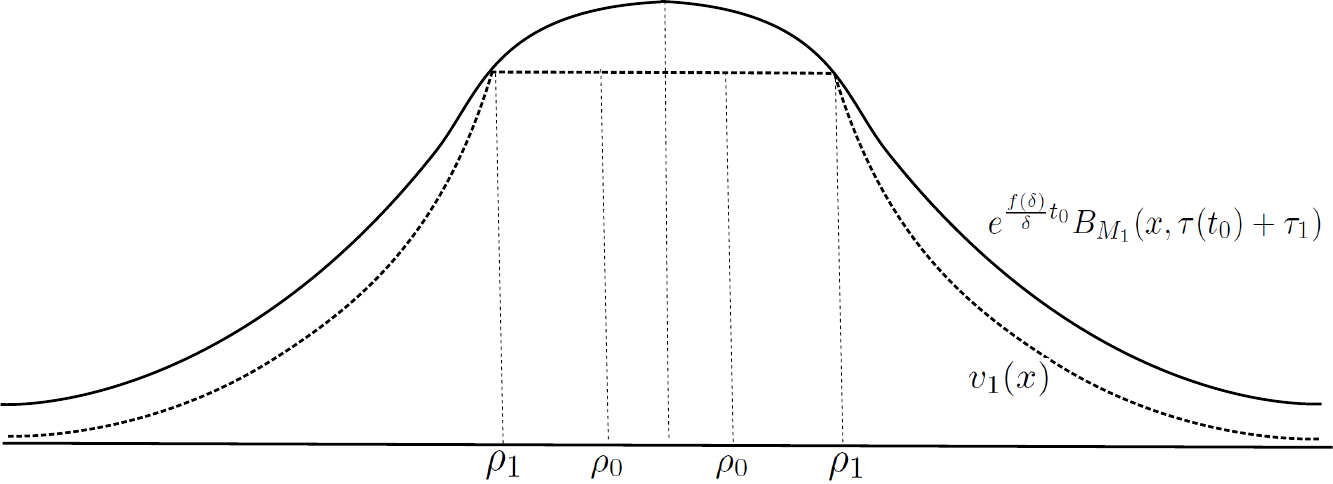}
  \caption{Step $IV$}
  \label{figv1}
\end{figure}
\noindent The proof is  complete for $m>1$ and $k=1$ (see Figures \ref{figv0} and \ref{figv1} for the construction of $v_1$).
\medskip

\noindent {VI. \sl The iteration.} We are now ready to address the next delicate step. Once we have proved that  $v(x,t_0) \geq v_1(x)$  for all $x\in \RN$, where $v_1$ is defined above, we apply the same proof and result to obtain
$$
v(x,2t_0) \geq (\text{solution of KPP with initial data }v_1(x))(t_0) \geq v_2(x),
$$
where $v_2(x)$ has the same construction as $v_0(x)$ and $v_1(x)$ but with parameters $\rho_2$ and $a_2$. Since $\rho_1> \rho_0$ the previous choice of $t_0$ is still valid to get to a similar conclusion. The argument continues for all $k=3,4,...$.

Let us check more closely the quantitative part of the iteration  in order to get an improvement. In the  process we keep $\epsilon$ fixed but me replace $\rho_0$ by $\rho_k$, $k\ge 1$, so that the formula \eqref{rho1} becomes
$$
\frac{\rho_{k+1}}{\rho_k}= L_k := (K_2/2K_1)^{1/(N+2s)}\, e^{\frac{f(\delta)}{(N+2s)\delta}t_0}
\left(1+( \tau_0/\tau_1(\rho_k))\right)^{2s\beta/(N+2s)}\,.
$$
As $k\to\infty$ we have $\rho_k\to \infty$, hence $\tau_1(\rho_k)\to \infty$, and the last quantity tends to
$$
L_\infty= (K_2/2K_1)^{1/(N+2s)}\, e^{\frac{f(\delta)}{(N+2s)\delta}t_0}\,.
$$
Finally, if we are given some $\sigma <\sigma_2=f'(0)/(N+2s) $ we can change the definition of $t_0$ so that we also have $L_\infty\ge e^{\sigma t_0}$.
 The conditions we put on $\delta$ and $t_0$ can be summarized in \eqref{delta} and \eqref{t0def1}, and they are independent on the parameter $\tau_k$, $\rho_k$ of the iteration.\normalcolor This ends the proof for  $m>1$.

\medskip

\noindent {\sc Case $m<1$.} The outline of the proof is similar to the case $m>1$. We explain the differences that appear.
The new time $\tau$ is introduced via
\begin{equation}\label{tauDef2}
\tau=\frac{1}{(1-m)f'(0)}\left[1-e^{-(1-m)f'(0)t}\right] \text{ if } m>1.
\end{equation}
Therefore, for each $t$ we have a new bounded time $\tau(t) \le \tau_{\infty}=1/((1-m)f'(0))$.  This property allows us to simplify the choice of $t_0$ as follows: condition \eqref{t0def1} is satisfied if
\begin{equation}\label{t0def2}
e^{\frac{f(\delta)}{\delta}t_0} \ge K_3 \left(1+( \tau_\infty/\tau_1(1))\right)^{N \beta},
\end{equation}
where $\tau_\infty/\tau(1) =\epsilon^{m-1}(c_2(1-m)f'(0))^{-1}$ is a constant independent of $\rho_0, \rho_1.$

Summing up: we take $\delta$ small enough such that
$$\frac{f(\delta)}{(N+2s)\delta} \geq \sigma$$
and $t_0$ such that
\begin{equation}\label{condt02}
  e^{\frac{f(\delta)}{\delta}t_0}\ge K_3 \left(1+( \tau_{\infty} /\tau_1(1))\right)^{N \beta},  \qquad (K_2/2K_1)^{1/(N+2s)} \, e^{\frac{f(\delta)}{(N+2s)\delta}t_0} \ge e^{\sigma t_0}.
\end{equation}
The rest is essentially the same.  \finprf

\medskip

%

%%%%%%%%%%%%%%%%%%%%%%%%%%%%%%%%%%%%%%%%%%%%%%%%%%%%%%%%%%

\begin{lemma} [ Expansion of uniform positivity for all times] \label{lemmaLowerBoundEps}
Let $N\geq 1$, $s\in (0,1)$, $m_1<m$, $f$ satisfying \eqref{propf} and $\sigma \in (0,\sigma_{2})$. Let $t_0>0$ from Lemma {\rm \ref{Lemmav0}}.
Then for every measurable initial datum $u_0$ with $0\leq u_0 \leq 1$, $u_0 \neq 0$, there exist $\epsilon \in (0,1)$ and $b>0$ (both depending on $u_0$) such that the solution $u$ of Problem \eqref{KPP} with initial datum $u(0,\cdot)=u_0$ satisfies
$$
u(x,t) \geq \epsilon \text{  for all }t \geq t_0 \text{ and }|x|\leq b e^{\sigma t}\,.
$$
\end{lemma}
\prf
Let $t_0$ defined in Lemma \ref{Lemmav0}. Then by Lemma \ref{lemmaTail} there exist $\epsilon>0$, $a_0>0$, $\rho_0>1$ such that $u(x,t)$ is bounded from below by a function $v_0$ with the long tail behavior at infinity
$$
u(x,t) \geq v_0(x):=\left\{
  \begin{array}{ll}
    a_0|x|^{-(N+2s)}, & \hbox{$|x|\geq \rho_0$,} \\
    \epsilon=a_0 \rho_0^{-(N+2s)}, & \hbox{$|x|\leq \rho_0$,}
  \end{array}
\right.$$
for all $x \in \RN, \ t \in [t_0,2t_0].$
In this way $v_0$ can be taken as the initial datum \eqref{v0def} in Lemma \ref{Lemmav0}.
We make $a_0$ smaller, if necessary, to have that $\epsilon=a_0 \rho_0^{-(N+2s)}\leq \epsilon_0$, where $\epsilon_0$ is given in Lemma \ref{Lemmav0}.

Therefore, by applying Lemma \ref{Lemmav0}, the solution $u$ will be raised an
$\epsilon$ at a large time $\tau_0+t_0$ and this holds true for all $\tau_0 \in [t_0,2t_0]$. More exactly, by \eqref{vxkt0}, for every $k=0,1,2,...$ one has
$$u(x,\tau_0 + kt_0) \geq \epsilon \text{ for all  }|x| \leq e^{\sigma k t_0} \rho_0, \ \tau_0 \in [t_0,2t_0]$$
which rewrites as
\begin{equation}\label{ineq1}
u( x,t) \geq \epsilon \text{ for all  }|x| \leq e^{\sigma k t_0}\rho_0, \quad t \in [(k+1)t_0, (k+2)t_0].
\end{equation}
But for $ t \in [(k+1)t_0, (k+2)t_0]$ we get $e^{\sigma k t_0}= e^{\sigma k t_0 - \sigma t}e^{\sigma t} \geq e^{-2\sigma t_0 }e^{\sigma t}$
and then \eqref{ineq1} implies, in particular, that
$$
u(x,t) \geq \epsilon, \text{ for all  }|x| \leq e^{-2\sigma  t_0}e^{\sigma t} \rho_0, \quad t \in [(k+1)t_0, (k+2)t_0].
$$
Since the union the intervals $[(k+1)t_0, (k+2)t_0]$ with $ k =0,1,2,....$ cover all $[t_0,\infty)$, we deduce that
$$
u(x,t) \geq \epsilon \quad \text{if} \quad t \geq t_0 \text{ and }|x|\leq \rho_0e^{-\sigma 2 t_0}e^{\sigma t}.
$$
The proof of the lemma follows by denoting $b=\rho_0 e^{-\sigma 2 t_0}$. \finprf

\smallskip

%%%%%%%%%%%%%%%%%%%%%%%%%%%%%%%%%%%%%%%%%%%%%%%%%%%%%

\begin{lemma}\label{lemmaLowerBoundEps2}
Let $N\geq 1$, $s\in (0,1)$, $f$ satisfying \eqref{propf}. Let $\sigma_{2}=\frac{f'(0)}{N+2s}$. Let $u$ be a solution of Problem \eqref{KPP} with initial datum $0\leq
u_0(\cdot)\leq 1$, $u_0 \neq 0$. Then for every $\sigma<\sigma_{2}$ there exist $\epsilon\in(0,1)$ and $\underline{t}>0$ such that
$$
u(x,t)\geq \epsilon \quad \text{ for all }t \geq \underline{t} \text{ and }|x|\leq e^{\sigma t}.
$$
\end{lemma}
\prf
We apply Lemma \ref{lemmaLowerBoundEps} with $\sigma$ replace by $\sigma' \in (\sigma, \sigma_2)$. Since $e^{\sigma t}\leq b e^{\sigma' t}$ for $t$ large, where
$b$ is the
constant in the statement of Lemma \ref{lemmaLowerBoundEps}, we deduce that
$$
u(x,t) \geq \epsilon \text{ for }t \geq \underline{t} \text{ and }|x|\leq e^{\sigma t}.
$$
\finprf

%%%%%%%%%%%%%%%%%%%%%%%%%%%%%%%%%%%%%%%%%%%%%%%%%%%%%%%%%%%%%%%%

\subsection{Case $m_c<m<m_1$}\label{subsectFirstRange}

In a similar way, we can prove the convergence to $1$ on the inner sets also in the range of parameters $m_c<m<m_1$.

\begin{proposition}\label{proConv1SecondCase}
Let $N\geq 1$, $s\in (0,1)$, $m_c < m<m_1$, $f$ satisfying \eqref{propf}. Let $\sigma_{1}=\frac{1-m}{2s}f'(0)$.
Let $u$ be a solution of Problem \eqref{KPP} with initial datum $0\leq u_0(\cdot)\leq 1$, $u_0 \neq 0$. Then for every $\sigma \in (0,\sigma_{1})$, $u(x,t)
\rightarrow 1$ uniformly on $\{|x|\leq e^{\sigma t}\}$ as $t \rightarrow \infty.$
\end{proposition}

\prf We argue in a similar way as in the case $m>m_1$ proved in Proposition \ref{proConv1}. The difference appears when obtaining the positivity on inner sets.
To this aim,  we start with nontrivial initial data $0\leq u_0\leq 1$ and we prove the analogue of Lemma \ref{lemmaLowerBoundEps}. The key ingredient is to use the
quantitative lower estimates for the solution $\underline{u}(x,t)$ Fractional Fast Diffusion Equation stated in Theorem \ref{ThGlobalLowerBounds1} to obtain an
estimate of the form
$$\underline{u}(x,t) \geq v_0(x), \quad \forall t \in [t_0,2t_0], \ x\in \RN,$$
where $v_0(x)$ is defined as
\begin{equation}\label{v0def2}
v_0(x)=\left\{
  \begin{array}{ll}
    a_0|x|^{-2s/(1-m)}, & \hbox{$|x|\geq \rho_0$,} \\
    \epsilon=a_0 \rho_0^{-2s/(1-m)}, & \hbox{$|x|\leq \rho_0$\,,}
  \end{array}
\right.\end{equation}
Afterwards, we can prove an analogue result to Lemma \ref{lemmaConv0} starting with initial data of the form \eqref{v0def2}.
Since the Barenblatt solution has a long tail decay of the form $|x|^{-2s/(1-m)}$, then we find $M_1>0$ and $\tau_1>0$ such that
$$
v_0(x) \geq B_{M_1}(x,\tau_1), \quad \forall x \in \RN.
$$
\finprf

%\subsection{Non-existence of traveling waves}\label{SubsecNonTravel}
%
%
%
%
%As a consequence of the exponential spread speed of the level sets we get that Problem \eqref{KPP} does not admit traveling waves solutions, i.e. solutions of the
%form
%$$u(x,t)=\varphi(x+te), \quad x \in \RN, t\in \RR,$$
%where $e$ is a unit vector in $\RN$. Thus
%$$
%\nabla \varphi \cdot e + (-\Delta)^s \varphi^m=f(\varphi), \quad \text{in }\RN.
%$$
%However, one can admit the existence of logarithmic traveling waves. ...... \textcolor{blue}{To finish!}

%%%%%%%%%%%%%%%%%%%%%%%%%%%%%%%%%%%%%%%%%%%%%%%%%%%%%%%%%%%%%%%%%%%%%%%%%%%%%%%%%

\section{The linear diffusion problem }\label{sectSelfSimLinProb}

We will need a number of facts about the linear diffusion equation for $0<s<1$,
\begin{equation}\label{LinProblem1}
  U_t + (-\Delta)^s U=0 \quad \text{for } x \in \mathbb{R}^N \text{ and }t>0.
\end{equation}
This problem has been studied, mainly by probabilists (\cite{Applebaum, Bertoin}), see also \cite{Valdinoc}, and many results are known. When considering initial
data $U_0 \in L^1(\ren)$, or more general,
\begin{equation}\label{LinProblemData}
 U(0,x)=U_0(x) \quad \text{for } x \in \mathbb{R}^N\,,
\end{equation}
the solution of Problem  \eqref{LinProblem1}-\eqref{LinProblemData} has the integral representation
$$
U(x,t)=\int_{\RN}K_s(x-z,t)U_0(z)dz\,,
$$
where the kernel $K_s$ has Fourier transform $\widehat{K}_s(\xi,t)=e^{-|\xi|^{2s}t}.$
If $s=1$, the function $K_1(x,t)$ is the Gaussian heat kernel.

\subsection{The fundamental solution. Further results on the asymptotics for large $|x|$}

We need some detailed information on the behaviour of the kernel $K_s(x,t)$ for $0<s<1$. In the particular case $s=1/2$, the kernel is explicit, given by the
formula
$$
K_{1/2}(x,t)=C_N t (|x|^2+t^2)^{-(N+1)/2}.
$$
In general, we know that the kernel $K_s(x,t)$ is the fundamental solution of Problem \eqref{LinProblem1}, that is $K_s(x,t)$ solves the problem with initial data
the Delta function
$$
\lim_{t \rightarrow 0}K_s(x,t)=\delta(x).
$$
It is known that the kernel $K_s$ has the form
$$K_s(x,t)=t^{-N/2s}f(|x|t^{-1/2s})$$
for some profile function, $f(r)$, that is positive and decreasing, and behaves at infinity like $f(r) \sim r^{-(N+2s)}$, cf. \cite{BlumenthalGetoor}.

\noindent We perform now a further analysis of the properties of the fundamental solution. Our aim is to prove the following result.

\begin{proposition}\label{proDerivKernel}
For every $s\in (0,1)$, the fundamental solution $K_s(x,t)$ of Problem \eqref{LinProblem1} is a increasing function in time
$$
\frac{d}{dt}K_s(x,t) \geq 0, \quad \text{for all large values of   }\ |x|/t^{1/2s}.
$$
\end{proposition}

 This property is known to be satisfied for the fundamental solution of various types of diffusion equations of evolution type: the Gaussian profile for the Heat
 Equation, the Barenblatt solution for the Fast Diffusion Equation.

The analysis of the derivative $\frac{d}{dt}K_s(x,t)$ involves not only the characterization of the profile $f$ for large $r$, but also a similar property for the
derivative $f'$. In fact, we will prove that $f(r)$ and $rf'(r)$ have the same behavior for large arguments. This is due to the power decay property of the
profile $f$.

 We recall that this property  is clearly true in the explicit case $s=1/2$ where $f(s)=(1+s^2)^{-(N+1)/2}$. But it is not true in the limit $s\to 1$, i.\,e., in
 the case of the Gaussian profile of the Heat Equation $G(s)=e^{-s^2/4}.$ Indeed, we can not obtain the same behavior for $G(s)$ and $s G'(s)$ since in this case
 the profile has an exponential expression.

\medskip

\noindent {\sl Proof of the proposition.} We recall that
\begin{equation}\label{scalingf}
K_s(x,t)=t^{-\frac{N}{2s}}f_{2s}(1,t^{-\frac{1}{2s}}|x|)
\end{equation}
(\cite{BlumenthalGetoor}), where $f_{2s}(1,x)$ is a continuous strictly positive function on $\RN$ of radial type, which is explicitly given by the expression
\begin{align*}
f_{2s}(1,x)&=\left[ (2\pi)^{N/2} |x|^{\frac{N}{2}-1} \right]^{-1} \int_0^{\infty} e^{-\omega^{2s}}\omega^{\frac{N}{2}}J_{\nu}(|x|\omega)d\omega \\[2mm]
&=\frac{1}{(2\pi)^{N/2} |x|^{N} } \int_0^{\infty} e^{-\left(\frac{\omega}{|x|}\right)^{2s}}\omega^{\frac{N}{2}}J_{\nu}(\omega)d\omega, \quad \nu=(N-2)/2 ,
\end{align*}
where $J_\mu$ denotes the Bessel function of first kind of order $\mu$. For simplicity, we denote $f(r)=f_{2s}(1,x), \ r=|x|$ since $f_{2s}(1,\cdot)$ is a radial
function:
\begin{equation}\label{f2s}
f(r)= \frac{1}{(2\pi)^{N/2}} \ r^{-N}  \int_0^{\infty} e^{-\left(\frac{\omega}{r}\right)^{2s}}\omega^{\frac{N}{2}}J_{\nu}(\omega)d\omega, \quad \nu=(N-2)/2.
\end{equation}

\noindent Next, we prove an intermediate result, concerning the behavior of the derivative $f'$.

%\medskip

%\medskip

\begin{lemma}\label{lemma1DerivFundSol}
Let $s\in (0,1)$ and let $f(r)=f_{2s}(1,x)$ be defined by \eqref{f2s}. Then
$$\lim_{r \rightarrow \infty} r^{N+2s} (Nf(r)+ r f'(r)) = - s^2 2^{2s+1} \frac{1}{\pi^{1+N/2}  } (\sin \pi s ) \Gamma(s)\Gamma\left(s+\frac{N}{2}\right). $$
In particular, we prove that $rf'(r) \sim -r^{-(N+2s)}$ for large $r$.
\end{lemma}
\prf
We compute the derivative with respect to $r$
%f'(r)= - (2\pi)^{-N/2} r^{-\mu-1} \int_0^{\infty} e^{-t^{2s}}t^{\frac{N}{2}} rt \ J_{\mu+1}(rt)dt, \quad \mu=\frac{N}{2}-1.
$$
f'(r)=  \frac{1}{(2\pi)^{N/2}  } \ r^{-N-1} \int_0^{\infty} \left( -N+2s\left(\frac{\omega}{r}\right)^{2s} \right)
 e^{-\left(\frac{\omega}{r}\right)^{2s}}\omega^{\frac{N}{2}}J_{\nu}(\omega)d\omega.$$
Therefore
\begin{align*}
r f'(r)&=  -Nf(r) + \frac{1}{(2\pi)^{N/2}  } \ r^{-N} \int_0^{\infty} 2s\left(\frac{\omega}{r}\right)^{2s} \,
e^{-\left(\frac{\omega}{r}\right)^{2s}}\omega^{\frac{N}{2}}J_{\nu}(\omega)d\omega =(I) +(II)\,,
\end{align*}
where (I)=$-Nf(r)$, and (II) is given by
 $$(II)=2s\frac{1}{(2\pi)^{N/2}  } \ r^{-(N+2s)} \int_0^{\infty}
 e^{-\left(\frac{\omega}{r}\right)^{2s}}\omega^{2s+\frac{N}{2}}J_{\nu}(\omega)d\omega\,.
$$
According to formula \eqref{Besel1}, we can write
$$\omega J_{\frac{N}{2}-1}(\omega) =NJ_{\frac{N}{2}}(\omega) - \omega J_{\frac{N}{2}+1}(\omega),$$ and therefore
\begin{align*}
 (II) &= 2Ns\frac{1}{(2\pi)^{N/2}  } \ r^{-(N+2s)} \int_0^{\infty} e^{-\left(\frac{\omega}{r}\right)^{2s}}\omega^{2s+\frac{N}{2}-1}
 J_{\frac{N}{2}}(\omega)  d\omega \\
  &- 2s\frac{1}{(2\pi)^{N/2}  } \ r^{-(N+2s)} \int_0^{\infty} e^{-\left(\frac{\omega}{r}\right)^{2s}}\omega^{2s+\frac{N}{2}}
J_{\frac{N}{2}+1}(\omega) d\omega
\end{align*}
Then, according to P\'olya (see Blumenthal \cite{BlumenthalGetoor})
$$\lim_{r \rightarrow \infty} \int_0^{\infty} e^{-\left(\frac{\omega}{r}\right)^{2s}}\omega^{2s+\frac{N}{2}-1}
 J_{\frac{N}{2}}(\omega)  d\omega = \frac{2}{\pi} \sin \pi s \int_0^{\infty} \omega^{2s+\frac{N}{2}-1}
 K_{\frac{N}{2}}(\omega)  d\omega.$$
and
$$\lim_{r \rightarrow \infty} \int_0^{\infty} e^{-\left(\frac{\omega}{r}\right)^{2s}}\omega^{2s+\frac{N}{2}}
J_{\frac{N}{2}+1}(\omega) d\omega = \frac{2}{\pi} \sin \pi s \int_0^{\infty}\omega^{2s+\frac{N}{2}}
 K_{\frac{N}{2}+1}(\omega)  d\omega.$$
Here the function $K_{\mu}$ are described in the paper of Erd\'elyi \cite{Erdelyi} (not to be confused with $K_s(x,t)$). Moreover (\cite{Erdelyi} page 51) we have
$$
L_1= \int_0^{\infty} \omega^{2s+\frac{N}{2}-1}  K_{\frac{N}{2}}(\omega)  d\omega=2^{2s+\frac{N}{2}-2}\Gamma\left(s+\frac{N}{2}\right)\Gamma(s).
$$
 $$
 L_2=\int_0^{\infty} \omega^{2s+\frac{N}{2}}  K_{\frac{N}{2}+1}(\omega)  d\omega=2^{2s+\frac{N}{2}-1}\Gamma\left(s+\frac{N}{2}+1\right)\Gamma(s).
$$
Therefore,
\begin{align*}\lim_{r\rightarrow \infty}& r^{N+2s} \left(r f'(r) + N f(r)\right)=-2s C_1(N,s),
\end{align*}
where
\begin{equation}\label{C1}
C_1(N,s):=s 2^{2s} \frac{1}{\pi^{1+N/2}  } (\sin \pi s ) \Gamma(s)\Gamma\left(s+\frac{N}{2}\right).
\end{equation}
If we write this result as
$$r^{N-1}\left(r f'(r) + N f(r)\right) \sim - 2sC_1(N,s)r^{-2s-1}$$
by integrating we obtain $r^{N}f(r) \sim C_1(N,s) r^{-2s}, $ that is
$$f(r) \sim C_1(N,s) r^{-(N+2s)},$$
which is exactly the result proved in \cite{BlumenthalGetoor}.
Moreover, we obtain that
$$\lim_{r\rightarrow \infty} r^{N+2s} r f'(r)= -(N +2s)C_1(N,s),$$
that is
$$r f'(r) \sim -  r^{-(N+2s)} \quad \text{for large }r.$$
\finprf

We complete the proof of Proposition \ref{proDerivKernel} on the behavior of the fundamental solution for large values of $\eta= |x|\,t^{-1/2s}.$

\prf
The Fundamental solution is given by
$$K_s(x,t)=t^{-\frac{N}{2s}}f(t^{-\frac{1}{2s}}|x|) .$$
We compute the derivative in the $t$ variable. According to the scaling formula \eqref{scalingf} we obtain
\begin{align*} \frac{d}{dt}K_{s}(x,t) &=  -\frac{N}{2s} t^{-\frac{N}{2s}-1}f(t^{-\frac{1}{2s}}|x|) -
\frac{1}{2s}t^{-\frac{N}{2s}-\frac{1}{2s}-1}|x|f'(t^{-\frac{1}{2s}}|x|) \\
&=-\frac{1}{2s} t^{-\frac{N}{2s}-1} \left[  Nf(\eta) + \eta f'(\eta) \right], \quad \eta =t^{-\frac{1}{2s}}|x|.
\end{align*}
By Lemma \ref{lemma1DerivFundSol} we know that
$$ Nf(\eta) + \eta f'(\eta)  \sim -2s C_1(N,s) \eta^{-(N+2s)}, \quad \text{for large }\eta,
$$
where $C_1(N,s)$ is a positive constant given  by formula \eqref{C1}. Therefore,
$$
\frac{d}{dt}K_{s}(x,t) \sim  t^{-\frac{N}{2s}-1} C_1(N,s) \eta^{-(N+2s)}= C_1(N,s)|x|^{-(N+2s)}, \quad \text{for large }\eta.
$$
\finprf

%%%%%%%%%%%%%%%%%%%%%%%%%%%%%%%%%%%%%%%%%%%%%%%%%%%%%%%%%%%%%%%%%%%%%%%%%%

\subsection{Self-similar solutions of the linear diffusion problem}

We study the existence, uniqueness and properties of self-similar solutions  of the form
\begin{equation}\label{selsim}
U(x,t)=t^{\alpha_1}F(t^{\beta_1}|x|)
\end{equation}
of the linear problem (the FPM Equation)
\begin{equation}\label{LinProblem}
  \left\{ \begin{array}{ll}
  U_t + (-\Delta)^s U=0 &\text{for } x \in \mathbb{R}^N \text{ and }t>0, \\
  U(0,x)=U_0(x)=C\,|x|^{\gamma}. &\text{for } x \in \mathbb{R}^N,
    \end{array}
    \right.
\end{equation}
where $C>0$, and $0<\gamma<2s$ is given. The constants $\alpha_1, \beta_1 \in \mathbb{R}$ will be determined such that $U(x,t)$ is a self-similar solution of
Problem \eqref{LinProblem}.

\smallskip

\noindent \textbf{Existence} of a solution $U$ to Problem \eqref{LinProblem} follows from  paper \cite{BV2012}, since the initial data $U_0(x)=|x|^{\gamma}$ with
$\gamma<2s$ belongs to a suitable weighted space $L^1(\RN,\varphi dx)$.

Let $\eta=t^{\beta_1}|x|$. Then,
$$U_t(x,t)=\alpha_1 t^{\alpha_1-1}F(\eta) + \beta_1 t^{\alpha_1-1}\eta F'(\eta), $$
$$(-\Delta)^s U(x,t)=t^{\alpha_1}(-\Delta)^s (F(t^{\beta_1}|x|))  = t^{\alpha_1}t^{2\beta_1 s}(-\Delta)^s F(\eta).$$
We obtain a first relation on the parameters: $\alpha_1-1=\alpha_1+2\beta_1 s$, and then $\beta_1=-\frac{1}{2s}$.

\smallskip

\noindent \textbf{Equation. }The profile $F$ satisfies the equation
$$
\alpha_1 F(\eta) + \beta_1 \eta F'(\eta) + (-\Delta)^s F(\eta)=0.
$$

\noindent\textbf{Self-similarity condition.} The equation is invariant under transformations of the form
$$
T_{\lambda}U(x,t)=\lambda^{-\alpha_1} U (\lambda^{-\beta_1}x, \lambda t).
$$
Therefore, $u=T_{\lambda}u.$ We apply this to the initial data
$$
T_{\lambda}U(x,0)=\lambda^{-\alpha_1} U (\lambda^{-\beta_1}x,0)=\lambda^{-\alpha_1-\beta_1 \gamma} |x|^{\gamma}
$$
and then $ \alpha_1=-\gamma \beta_1. $ We obtain the exact value of the similarity exponents
\begin{equation}\label{AlphaBeta}
\alpha_1=\frac{\gamma}{2s}, \quad \beta_1=-\frac{1}{2s}.
\end{equation}
Notice that $\alpha_1>0$ and $\beta_1<0$. As a solution of the linear problem \eqref{LinProblem}, $U(x,t)$ can be computed as a convolution with the kernel
$K(\cdot,t)$
$$
U(x,t)= (K(\cdot,t) \star U_0)(x)=\int_{\RN}K(y,t)U_0(x-y) dy.
$$
Since the initial data is a radial function $U_0(x)=|x|^{\gamma}$, then by the properties of the kernel $K$, $U$ will also be a radial function, and therefore the
profile $F$ is radial.

\begin{lemma}[\textbf{Properties of the profile}]\label{lemmaSelfSim}
 The profile $F$ is monotone non-decreasing and it satisfies $\eta F' \leq c_2 F$, for all $\eta \geq 0.$
\end{lemma}
\prf

\noindent \textbf{I. Monotonicity property. }In order to prove the positivity of $F$ we will make use of the Alexandrov Symmetry Principle and we prove that
$U(x,t)$ is radially increasing in the space variable $x \in \RN.$

We start with increasing radial initial data $U_0(x)=|x|^{\gamma}$. We approximate $U_0$ with a sequence of radially symmetric and bounded functions $U_{0n} \in
L^1(\RN)$ such that $U_{0n}(r)\to C\, n^{\gamma}$ as $r\to\infty$ and $v_{0n}(r)= C\,n^{\gamma}- U_{0n}(r)\in L^1(\RN)$. Let $v_n$ the solution of Problem
\eqref{LinProblem} with initial datum $v_{0n}$. We may apply the Alexandrov Symmetry Principle (that we explain  in detail below) to $v_n$ to conclude that it is
radially symmetric and decreasing w.r.t. the space variable. We then put $U_n(x,t)=C\,n^{\gamma}-v_n(x,t)$, which is  radially symmetric and increasing, and
solves \eqref{LinProblem} with initial datum $U_{0n}$. We pass now to the limit $n\to\infty$ to get the same conclusion for $U$.

\smallskip

\textbf{Applying the Alexandrov Symmetry Principle. }We fix two points $x$ and $x'$ in $\RN$ such that $|x|<|x'|$. Let $H$ denote the hyperplane perpendicular on
the line $xx'$.
Let $\Omega_1$ and $\Omega_2$ be the two sets delimited by the hyperplane $H$ such that the origin is contained in $\Omega_1$. Let $\Pi$ the symmetry
 with respect to $H$ that maps $\Omega_1$ into $\Omega_2$. Clearly, $\Pi(x)=x'$, $x \in \Omega_1$.
Then one can prove that for every $y \in \Omega_1$ $|y|<|y'|$, where $y'=\Pi (y)$. Since $v_{0n}$ is radially decreasing, we get that $v_{0n}(y) \geq
v_{0n}(\Pi(y))$, for all $y \in \Omega_1$. By applying the Alexandrov Symmetry Principle stated in Theorem \ref{AlexPrincThm} we obtain that $v_n(x) \geq
v_n(x')$. The arguments we used can be done for every pair of points $|x|<|x'|$, therefore
 $v_n$ is radially increasing.

\medskip

\noindent \textbf{II. Decay at infinity. }A formal computation starting from the initial data  \ $U(x,t) \rightarrow |x|^{\gamma}$ as $t\rightarrow 0$ gives us
that
$ \eta^{-\gamma}F(\eta)  \rightarrow 1$ as $ \eta \rightarrow \infty$. Therefore
$$
F(\eta) \sim \eta ^{\gamma} \quad \text{ for large  }\eta.
$$
This characterization of the profile $F$ gives us the following spatial decay for $U(x,t)$ for large times
$$
C_2 |x|^{\gamma} \le U(x,t) \leq C_1 |x|^{\gamma}, \quad \text{for large } |x|t^{-1/2s}.
$$
Moreover, we will prove the following relation between $F'$ and $F$:
$$
|\gamma F(\eta) - \eta F'(\eta)| \le \eta^{\gamma}, \quad \text{for large }\eta >0 .
$$
As a consequence we can characterize the derivative $U_t$:
$$
U_t(x,t)=t^{\alpha_1-1}\frac{1}{2s} \left( \gamma F(\eta) - \eta F'(\eta) \right), \quad \eta=t^{-1/2s}|x|.
$$
$$
 U_t(x,t) \sim t^{-1}|x|^{\gamma} \quad \text{for large values of }\ t^{-1/2s}|x|.
$$
The first step will be to obtain a formula for the profile $F(\eta)$. Therefore
\begin{align*}U(x,t)&=K_s(x,t) \star U_0(x)= t^{-\frac{N}{2s}} \int_{\RN}f(t^{-\frac{1}{2s}}|x -y|) |y|^{\gamma}dy , \quad z=t^{-\frac{1}{2s}}y\\
&=t^{\frac{\gamma}{2s}}\int_{\RN}f(t^{-\frac{1}{2s}}x -z) |z|^{\gamma}dz.
\end{align*}
Since $U(x,t)$ has the self similar form \eqref{selsim} then
$$F(t^{-\frac{1}{2s}}x)=\int_{\RN}f(t^{-\frac{1}{2s}}x -z) |z|^{\gamma}dz = (f\star U_0)(t^{-\frac{1}{2s}}x) , \quad \forall x \in \RN, t>0,$$
that is
$$
F(\eta)=\int_{\RN}f(\eta-z) |z|^{\gamma}dz , \quad \forall \eta>0.
$$
Let us continue using the notations
$$
F(|\eta|)=F(\eta),\quad f(|\eta|)=f(\eta).
$$
We fix $\eta \in \RN.$ Let $|\eta|=\bar{\eta}$ and $\eta=\bar{\eta} e$ for a vector $e \in \RN$ with $|e|=1.$
Then
$$
F(\bar\eta)=\int_{\RN}f(|z|) |\eta-z|^{\gamma}dz=\bar\eta^{N+\gamma} \int_{\RN}f(|\bar \eta y |) |e-y|^{\gamma}dy, \quad z=\bar\eta y.
$$
We differentiate in $\bar{\eta}$
$$
 F'(\bar\eta)= \bar\eta^{N+\gamma-1} \int_{\RN} \left[ (N+\gamma) f(|\bar \eta y |) +  \bar \eta y   f'(|\bar \eta y |)      \right]|e-y|^{\gamma}dy.
 $$
Therefore
$$ \bar\eta F'(\bar\eta) - \gamma F(\bar\eta)  = \bar\eta^{N+\gamma} \int_{\RN} \left[ N f(|\bar \eta y |) +  \bar \eta y   f'(|\bar \eta y |)
\right]|e-y|^{\gamma}dy, \quad z=\bar\eta y$$
$$\bar\eta F'(\bar\eta) - \gamma F(\bar\eta)  = \bar\eta^{\gamma}\int_{\RN} \left[ N f(|z|) +  z  f'(|z |) \right] \left| e- \frac{z}
{\bar{\eta}}\right|^{\gamma}d z.$$
We know that
$N f(r) +  r  f'(r) \sim -C_1 r^{-(N+2s)}$ for large $r$. Since we deal with a convolution we will use the information only in the sense of modulus.  We fix $R>0$
such that
$$
C_1 r^{-(N+2s)} \leq |N f(r) +  z  f'(r) | \leq C_2 r^{-(N+2s)} , \quad \forall r\geq R.
$$
Then
$$
\int_{\RN} \left[ N f(|z|) +  z  f'(|z |) \right] \left| e- \frac{z} {\bar{\eta}}\right|^{\gamma}d z = I +II,
$$
where
\begin{align*}I= \int_{|z|\leq R}\left[ N f(|z|) +  z  f'(|z |) \right] \left| e- \frac{z} {\bar{\eta}}\right|^{\gamma}d z \leq C \int_{|z|\leq R}\left| e-
\frac{z} {\bar{\eta}}\right|^{\gamma}d z \leq C \left( 1+ \frac{R} {\bar{\eta}}\right)^{\gamma} R^N.
\end{align*}
The second term is estimated as follows (notice that $R$ can be taken large enough to ensure that $R>\bar{\eta}$ since the value $\bar{\eta}$ was fixed at the
beginning).
\begin{align*}II&= \int_{|z|\geq R}\left[ N f(|z|) +  z  f'(|z |) \right] \left| e- \frac{z} {\bar{\eta}}\right|^{\gamma}d z
\leq \int_{|z|\geq R} |z|^{-(N+2s)}  \left| e- \frac{z} {\bar{\eta}}\right|^{\gamma}d z \\
&\leq \int_{|z|\geq R} |z|^{-(N+2s)}  \left( 1 + \frac{|z|} {\bar{\eta}}\right)^{\gamma}
\leq 2 \int_{|z|\geq R} |z|^{-(N+2s)}  \left( \frac{|z|} {\bar{\eta}}\right)^{\gamma}, \\
&= 2 \bar\eta^{-\gamma} \int_{|z|\geq R} |z|^{-(N+2s)+\gamma} dz,  \quad \text{we know  } \gamma <2s,\\
&= 2 \frac{1}{\bar\eta^{\gamma}} \ \frac{1}{2s-\gamma} \ R^{\gamma-2s}.
\end{align*}
We conclude that
$$I +II  \leq C_1  + C_2 \frac{1}{\bar\eta^{\gamma}}.$$
Now, recall that $\bar\eta F'(\bar\eta) - \gamma F(\bar\eta)  =  \bar\eta^{\gamma}(I+II)$. Therefore we have proved that
$$  \limsup_{\bar\eta\rightarrow \infty} \bar\eta^{-\gamma} \left| \bar\eta F'(\bar\eta) - \gamma F(\bar\eta)  \right|  \leq C_1  + C_2
\frac{1}{\bar\eta^{\gamma}} \rightarrow C_1,
$$
which means the limit is finite.  Therefore,
$$
\left| \bar\eta F'(\bar\eta) - \gamma F(\bar\eta)  \right|  \le \bar\eta^{\gamma} \quad \text{for large }\bar\eta.
$$
\finprf

%%%%%%%%%%%%%%%%%%%%%%%%%%%%%%%%%%%%%%%%%%%%%%%%%%%%%%%%%%%%%%%%%%%%%%%%%%%%%%%%%%%%%%%%%%%%%%%

\section{The Reaction Problem}

As a further evidence of the influence of the tail of the data on the propagation rate, we consider the purely reactive problem (no diffusion)
\begin{equation}\label{reactionProblem}
u_t= f(u), \quad x \in \RN, \ t >0\,,
\end{equation}
with initial datum $u_0$  and  $f(u) \sim u(1-u)\sim f'(0)u.$
It is easy to see that when we simplify $f(u)$ to $f'(0)u=au$, the  exact solution is
$$
u(x,t)=u_0 e^{f'(0)t}.
$$
Let us examine the level sets in two particular cases.

\smallskip

\noindent\textbf{Exponential decay.} By considering initial datum of the form $u_0(x) \sim e^{-x^2}$ for large $|x|$, then the solution $u(x,t)$ satisfies a
similar behavior
$$
u(x,t) \sim e^{-(x^2-at)} \text{ for large } x.
$$
The level sets $u(x,t)=$ constant are characterized by  $x=\sqrt{at+c}$.

\noindent\textbf{Power decay.} By considering initial datum of the form $u_0(x) \sim |x|^{-(N+2s)}$ for large $|x|$, then the solution $u(x,t)$ is such that
$$
u(x,t) \sim e^{at}|x|^{-(N+2s)}.
$$
The level sets $u(x,t)=$ constant are characterized by $|x|\sim e^{\frac{a}{N+2s}t}$.
\medskip

\noindent\textbf{Conclusion: the influence of fractional diffusion}: For $|x|$ large, the solution of the reaction-diffusion Problem \eqref{KPP} behaves like the
solution of Problem \eqref{reactionProblem}, that is, the non-diffusion case. The fractional diffusion term $(-\Delta)^s u$ does no change the basic behaviour of
the solution for large $|x|$. This fact has been also observed by King and McCabe in \cite{KingMcCabe} in the fast diffusion case with the standard Laplace
operator.

%%%%%%%%%%%%%%%%%%%%%%%%%%%%%%%%%%%%%%%%%%%%%%%%

\section{Appendix}

\subsection{Concept of solution to Problem \eqref{KPP}}

 According to \cite{DPQRV2}  there exists a unique mild solution of Problem \eqref{genfractPME} corresponding to the initial datum $u_0 \in L^1(\RN)$, $0\leq u_0
 \leq 1$, constructed by means of the tools of semigroup theory. Moreover, such $u$ is in fact a strong solution of the equation. In the case $m>1$, the
 $C^{\alpha}$ regularity of the solution follows from \cite{AthanCaff}, and this has been extended to $m<1$ up to the extinction time (if there is one).
 Quantitative estimates of positivity of the solution for any $m>0$ corresponding to non-negative data have been proved in \cite{BV2012}. Recently, $C^{1,\alpha}$
 regularity of strong solutions was proved in \cite{DPQRV3}.

As a consequence one obtains by rather standard methods the existence, uniqueness and regularity properties of the solution to Problem \eqref{KPP} corresponding
to the initial datum $u_0 \in L^1(\RN)$, $0\leq u_0 \leq 1$. In order to prove the existence of a solution of the problem $u_t+L_s u^m=c u$, the idea is to prove that the map $u_0 \mapsto v=e^{-ct}u$ is a $m$-$\omega$-accretive operator. Standard properties, like the maximum principle hold also in our setting.

 A more detailed analysis of these properties is beyond the purpose of this work.

%%%%%%%%%%%%%%%%%%%%%%%%%%%%%%%%%%%%%%%%%%%%%%%%%%%%%%%%%%%%%%%%%%%%%

\subsection{A version of the Maximum Principle}

We need an interesting version of Maximum Principle  proved by Cabr\'{e} and Roquejoffre in \cite{CabreRoquejoffreArxiv}, Lemma 2.9, suitable for comparisons in
which fractional laplacian operators are involved.

\begin{lemma}\label{MaxPrinc}
Let $N\geq 1$, $s\in(0,1)$, $0\leq \gamma < 2s$. Let $v\in C^1([0,\infty);X_{\gamma})$ satisfy $v(\cdot,t)\in D_{\gamma}(A)$ for all $t>0$. Let $r:[0,\infty)$ be
a continuous function and define
$$\Omega_I=\{ (x,t)\in (0,\infty)\times \RN: |x|<r(t) \}. $$
Assume in addition:

  $(H1)$ \,  $v(0,\cdot)\leq 0$ in $\RN$.

  $(H2)$ \,  $v\le 0$ \ in \ $((0,\infty)\times \RN) \setminus \Omega_I.$

  $(H3)$  \, $a(x,t) v_t + L_s v \leq b\, v$ \ in \ $\Omega_I$.

Then $v\le 0$ in $(0,\infty)\times \RN$.
\end{lemma}

Although the equation we have is different, the proof as in \cite{CabreRoquejoffreArxiv} still works (with inessential modifications).

%%%%%%%%%%%%%%%%%%%%%%%%%%%%%%%%%%%%%%%%%%%%%%%%%%%%%%%%%%%%%%%%%%%%%%%%%%%%%%%%%

\subsection{Alexandrov Reflection Principle}\label{AlexPrincSect}

We recall the version of Alexandrov's symmetry principle that holds in the case of the nonlinear parabolic problem
\begin{equation}\label{fracPME}
u_t=L_s u^m,\quad u(0,x)=u_0(x),
\end{equation}
posed in $\RN$, with $L_s=(-\Delta)^s$, $m>0$, $u_0 \in L^1(\RN)$. Let us take a hyperplane $H$ that divides $\RN$ into two half-spaces $\Omega_1$ and $\Omega_2$
and consider the symmetry $\Pi$ with respect to $H$ that maps $\Omega_1$ into $\Omega_2$. The following result is proved as Theorem 15.2 in
\cite{VazquezBarenblattFractPME}:

\begin{theorem}\label{AlexPrincThm}
Let $u$ be the unique solution of Problem \eqref{fracPME} with initial data $u_0$. Under the assumption that
$$u_0(x) \geq u_0(\Pi(x)) \quad  \text{ in }\Omega_1$$
we have that for all $t>0$
$$u(x,t) \geq u(\Pi(x),t) \quad  \text{ for }x \in \Omega_1.$$
\end{theorem}

%%%%%%%%%%%%%%%%%%%%%%%%%%%%%%%%%%%%%%%%%%%%%%%%

\subsection{Bessel functions of first kind}
The Bessel function $J_\mu$ of first kind can be introduced through a series expansion, cf. \cite{Abramowitz},
$$
J_{\mu}(z)=\sum_{k=0}^{\infty} \frac{(-1)^k}{k! \ \Gamma(k+\mu+1)} \left( \frac{z}{2} \right)^{2k+\mu}.
$$
We mention the following recurrence formulas:
$$
J'_{\mu}(z)=\frac{1}{2} \left(  J_{\mu-1}(z) - J_{\mu+1}(z)\right), \quad \text{for  } \mu \ne 0.
$$
$$J'_{0}(z)=-J_{1}(z),$$
\begin{equation}\label{Besel1}
J_{\mu}(z)=\frac{z}{2\mu} \left(  J_{\mu-1}(z) + J_{\mu+1}(z)\right).
\end{equation}
\begin{equation}\label{Besel2}
\int_0^\infty K_a(t)t^{b-1}dt=2^{b-2}\Gamma\left( \frac{b+a}{2}\right)\Gamma\left( \frac{b-a}{2}\right), \quad \Re(b\pm a)>0.
\end{equation}

\

\section*{Comments and Open problems}

 \noindent$\bullet$ There are critical values of the speed $\sigma$ which we do not cover in this work: $\sigma_1$ for $m_c<m<m_1$; $\sigma_2$ for $m_1<m\le1$; respectively, $(\sigma_2,\sigma_3)$ for $m>1$.
The analysis of those cases leads to long new developments.

\noindent$\bullet$ Is there a definite profile function that represents up to translation the shape
of the solution in the region where it varies in a marked way to join the level $u=1$ to the level $u=0$? Maybe for $s=1/2$ this question is easier.

\noindent$\bullet$ For reasons of length and novelty, the case $m<m_1$ is not studied. For the corresponding fractional fast diffusion equation there appears the phenomenon of extinction in finite
time. King and McCabe in \cite{KingMcCabe} give an idea on the asymptotics in this range of parameters.

\noindent$\bullet$ A detailed numerical treatment of these problems for the case $m \ne 1$ is needed, see in this respect \cite{FdTeso}.

\noindent$\bullet$ There are other interesting directions in this class of problems. Thus, in a recent paper \cite{CabreCoulonRoquejoffreCRM}, the authors investigate the model
$$
u_{t}(x,t) + A u(x,t)=\mu(x)u-u^2, \quad x \in \RN, \ t>0,
$$
where the function $\mu$ is supposed periodic in each spatial variable $x_i$ and satisfy $0 < \min \mu \leq \mu(x)$.

\section*{Acknowledgments} Both authors have been supported by the Spanish Project MTM2011-
24696.

%\bibliographystyle{plain}
%\bibliography{biblioKPP}

%\begin{thebibliography}{10}
%\small

\end{document}